\documentclass{amsart}\usepackage{amsfonts,amssymb,amsmath,bbm,curves}
\usepackage{allrunes}
\def\C{\Bbb{C}}
\def\k{\mathbbm{k}}
\def\N{\Bbb{N}}
\def\P{\Bbb{P}}\def\R{\Bbb{R}}\def\Z{\Bbb{Z}}

\def\li{\ \\ $\bullet$ }

\def\suml{\sum\limits}

\def\cupl{\mathop\cup\limits}\def\prodl{\mathop\prod\limits}

\newcommand{\quotient}[2]{{\left.\raisebox{0.4ex}{$#1$}\!\!\middle/\!\!\raisebox{-0.4ex}{$#2$}\right.}}

\def\tC{\tilde{C}}\def\tE{\tilde{E}}
\def\tp{\tilde{p}}

\def\al{\alpha}\def\ga{\gamma}\def\Ga{\Gamma}\def\be{\beta}\def\de{\delta}\def\De{\Delta}

\def\si{\sigma}

\def\cG{\mathop\mathcal G}

\def\u0{\underline{0}}

\def\one{{1\hspace{-0.1cm}\rm I}}

\def\ber{\begin{array}{l}}\def\eer{\end{array}}
\def\bpm{\begin{pmatrix}}\def\epm{\end{pmatrix}}
\def\beqm{\begin{multline}}
\def\eeqm{\end{multline}}
\def\bM{\begin{matrix}}\def\eM{\end{matrix}}
\def\bee{\begin{enumerate}}\def\eee{\end{enumerate}}
\def\bet{\begin{tabular}{cccccccc}}\def\eet{\end{tabular}}

\def\sset{\subset}\def\smin{\setminus}
\def\cp{cut\&paste}
\def\omp{ordinary multiple point}
\def\Nnd{Newton-non-degenerate}\def\ND{Newton diagram}

\def\bull{\vrule height .9ex width .9ex depth -.1ex }

\newcommand{\beq}{\begin{equation}}\newcommand{\eeq}{\end{equation}}
\newtheorem{Lemma}{Lemma}[section]\newcommand{\bel}{\begin{Lemma}}\newcommand{\eel}{\end{Lemma}}
\newtheorem{Example}[Lemma]{Example}\newcommand{\bex}{\begin{Example}\rm}\newcommand{\eex}{\end{Example}}
\newtheorem{Proposition}[Lemma]{Proposition}\newcommand{\bprop}{\begin{Proposition}}\newcommand{\eprop}{\end{Proposition}}
\newtheorem{Definition-Proposition}[Lemma]{Definition-Proposition}

\def\bpr{~\\{\em Proof.\ }}\newcommand{\epr}{$\bull$\\}
\newtheorem{Theorem}[Lemma]{Theorem}\newcommand{\bthe}{\begin{Theorem}}\newcommand{\ethe}{\end{Theorem}}
\newtheorem{Definition}[Lemma]{Definition}\newcommand{\bed}{\begin{Definition}}\newcommand{\eed}{\end{Definition}}
\newtheorem{Remark}[Lemma]{Remark}\newcommand{\beR}{\begin{Remark}\rm}\newcommand{\eeR}{\end{Remark}}
\newtheorem{Corollary}[Lemma]{Corollary}\newcommand{\bcor}{\begin{Corollary}}\newcommand{\ecor}{\end{Corollary}}
\newtheorem{Conjecture}[Lemma]{Conjecture}\newcommand{\bconj}{\begin{Conjecture}}\newcommand{\econj}{\end{Conjecture}}

\newcommand{\bin}[2]{\binom{#1}{#2}}

\newcommand{\freeh}[3]{\setcounter{xx}{#1}\addtocounter{xx}{#3}\curve(#1,#2,\value{xx},#2)
\setcounter{yy}{#2}\addtocounter{yy}{3}\curve(#1,\value{yy},\value{xx},\value{yy})
\setcounter{yy}{#2}\addtocounter{yy}{-3}\curve(#1,\value{yy},\value{xx},\value{yy})}

\newcommand{\freev}[3]{\setcounter{yy}{#2}\addtocounter{yy}{#3}\curve(#1,#2,#1,\value{yy})
\setcounter{xx}{#1}\addtocounter{xx}{3}\curve(\value{xx},#2,\value{xx},\value{yy})
\setcounter{xx}{#1}\addtocounter{xx}{-3}\curve(\value{xx},#2,\value{xx},\value{yy})}

\newcommand{\claw}[3]{\setcounter{yy}{#2}\addtocounter{yy}{-#3}\curve(#1,#2,#1,\value{yy})
\setcounter{xx}{#1}\addtocounter{xx}{5}\curve(#1,#2,\value{xx},\value{yy})
\setcounter{xx}{#1}\addtocounter{xx}{-5}\curve(#1,#2,\value{xx},\value{yy})
}

\def\blackboxh{\rule{15pt}{5pt}}

\def\isom{\xrightarrow{\sim}}
\def\Types{\ensuremath{\mathfrak{Types}_{(\C^2,0)}}}
\def\Spec{\ensuremath{\mathfrak{Spec}_{(\C^2,0)}}}

\newcommand{\tinyM}{\scriptstyle}
\newcommand{\tinyT}{\scriptsize}

\title[]{R\MakeLowercase{ecombination formulae for the spectrum of curve singularities and some applications}}
\author[]{D\MakeLowercase{mitry} K\MakeLowercase{erner}}
\address{Department of Mathematics, Ben Gurion University of the Negev, P.O.B. 653, Be'er Sheva 84105, Israel.}
\email{dmitry.kerner@gmail.com}
\date{\today}
\thanks{Most of the results were obtained in Oberwolfach during my stay as the Leibnitz Fellow (2011).
 Further, I was partially supported by the grant FP7-People-MCA-CIG, 334347.
 \\We appreciate the computational power of "Singular".
 It is inspiring to peer at the ancient runes.
 {\textarc{a \a b c d \d e f F g h \h i \i j k l m n \ing o p \p q r \r 
 } }
}

\begin{document}\newcounter{xx}\newcounter{yy}
\begin{abstract}
We obtain some recombination formulae for the spectra of (complex, reduced) plane curve singularities. As an application we prove: a generalization of Durfee's bound;  a generalization of Givental's bound; the multiplicity of the curve singularity is determined by its spectrum; for many curve singularities all the multiplicities of exceptional divisors of the resolution are determined by the spectrum; etc.
\end{abstract}
\maketitle
\setcounter{secnumdepth}{6} \setcounter{tocdepth}{3}

\section{Introduction}
Consider a (complex, reduced) plane curve singularity  $(C,0)\sset(\C^2,0)$.
The singularity spectrum is a strong topological invariant of $(C,0)$, \cite{Steenbrink1976}.
 It is a collection of rational numbers, $-1<\xi_1\le\cdots\le\xi_\mu<1$, here $\mu$ is the Milnor number of $(C,0)$.
 We use the additive notation: $\suml^\mu_{i=1}t^{\xi_i}$.

This work has originated from the following "experimental observation".
The formulae expressing $Spec(C,0)$ through the spectra of the branches of $(C,0)$ are cumbersome.
 But a slightly modified question has a neat answer.
Consider the
(reduced) tangent cone $T_{(C,0)}=\{l_1,..,l_r\}$. Accordingly, consider the tangential decomposition
$(C,0)=\cupl^r_{\al=1}(C_\al,0)$. Here $(C_\al,0)$ can be further reducible, but its tangent cone consists of one (multiple) line and the tangent cones are distinct for different parts, $T_{(C_\al,0)}\cap T_{(C_\be,0)}=\{0\}$.
For
each $\al$ the "directional approximation" of $(C,0)$ is defined as follows. Let $(D,0)\sset(\C^2,0)$ be an ordinary multiple point
of multiplicity $mult(C,0)-mult(C_\al,0)$, i.e. a collection of smooth pairwise non-tangent branches.
Assume $T_{(C_\al,0)}\cap T_{(D,0)}=\{0\}$. Then the approximation of $(C,0)$ in the direction $\{l_\al\}$
is $(C^{(\al)},0):=(C_\al\cup D,0)$. This germ is not unique, but its embedded topological singularity type is well defined.

Using these notions the "experimental observation" on the spectrum can be stated as:
\beq
Spec(C,0)=\sum_{\al=1}^r Spec(C^{(\al)},0)-(r-1)Spec(x^m-y^m),\quad m=mult(C,0)
\eeq
This relation might look unexpected, it becomes more transparent if one compares the first blowup of both sides
(schematically):
\beq\label{Eq.Spectral.Relation.Simplest.Case}
\begin{picture}(0,0)(200,20)\linethickness{0.8mm}\curve(0,0,60,0)
\thinlines\curve(5,20,9,14,10,0)\curve(15,20,12,14,11,0)
\curve(26,10,30,0,26,-10)\curve(34,10,30,0,34,-10)\curve(38,10,35,5,30,0,25,-5,22,-10)
\curve(45,20,49,14,50,0)\curve(55,20,51,14,50,0)\curve(45,5,50,0,55,-5)
\put(70,-2){$\equiv$}
\end{picture}
\begin{picture}(0,0)(110,20)\linethickness{0.8mm}\curve(0,0,60,0)
\thinlines\curve(5,20,9,14,10,0)\curve(15,20,12,14,11,0)
\freev{30}{-5}{10}\freev{50}{-5}{10}\put(70,-2){+}
\end{picture}
\begin{picture}(0,0)(20,20)\linethickness{0.8mm}\curve(0,0,60,0)
\thinlines
\curve(26,10,30,0,26,-10)\curve(34,10,30,0,34,-10)\curve(38,10,35,5,30,0,25,-5,22,-10)
\freev{50}{-5}{10}\freev{10}{-5}{10}
\put(70,-2){+}
\end{picture}
\begin{picture}(0,0)(-70,20)\linethickness{0.8mm}\curve(0,0,60,0)
\thinlines\freev{10}{-5}{10}\freev{30}{-5}{10}
\curve(45,20,49,14,50,0)\curve(55,20,51,14,50,0)\curve(45,5,50,0,55,-5)
\put(65,-2){$-$}\put(75,-2){$(r-1)$}
\end{picture}
\begin{picture}(0,0)(-180,20)\linethickness{0.8mm}\curve(0,0,60,0)
\thinlines\freev{10}{-5}{10}\freev{30}{-5}{10}\freev{50}{-5}{10}
\end{picture}
\eeq
\\\\\\
Here $\equiv$ denotes the equality of classes in an abelian group defined later.

\

The relation above is just the simplest among numerous relations of singularity spectrum. We call them the "recombination" or "\cp" formulas.
 In this paper we give a method to construct many relations of this type.
Given some singularity types, $\{C_i,0\}$ , and some integers, $\{a_i\}$, we write $\sum a_i(C_i,0)\equiv0$ if
 $\sum a_iSpec_{(C_i,0)}=0$. Usually we write this in terms of (partial) resolution, e.g. as in equation (\ref{Eq.Spectral.Relation.Simplest.Case}).

\

Consider the abelian group $\Spec$ of the linear combinations $\suml_i a_i Spec(C_i,0)$, generated by all the possible (reduced) plane curve singularities. Using recombination formulae, as above, we prove that $\Spec$ is freely generated by the very restricted class of singularities (we call them "basic types"),  $C_{(p,q)}:=\{(x^p+y^p)(y^q-x^{2q})=0\}\sset(\C^2,0)$. For $q=0$ or $q=1$ these are the ordinary multiple points, for $q>1$ there are smooth tangent branches as well.

Further, we introduce the natural product, so that $\Spec$ becomes a (commutative, associative, unital) ring, generated by $C_{(1,0)}$ and $\{C_{(0,q)}\}_{q\ge2}$.

\

A natural question is whether similar formulae hold for other (stronger) invariants. In section \ref{Sec.Spectral.Pairs.Non.additive} we show that even the simplest formula of equation (\ref{Eq.Spectral.Relation.Simplest.Case}) does not hold for the spectral pairs, neither for the signatures of the Seifert form. Thus the spectrum of plane curve singularity seems to be on the border between the "simple" invariants (easily computable and controlled) and very strong invariants.

The recombination formulae suggest many applications, we consider briefly few of them. Many proofs are missing, more detail will be given in the subsequent version of the paper.
\li An approach to the classification of surface singularities $\{z^2+f(x,y)\}$ with small values of $\mu_++\mu_0$, section \ref{Sec.Surfaces.Sings.low.muPlusZero}.
\li A generalization of Durfee bound for isolated hypersurface singularities $(X,0)\sset(\C^n,0)$, section \ref{Sec.Durfee.Bound.Generalized}.
For any $0\le\al<1$ and any $n\ge2$ there exists a constant $C_n$ (depending on $n$ only) such that
$\frac{\mu(X,0)}{n!}\ge\frac{|Spec(X,0)\cap(-1,-\al]|+C_n}{(1-\al)^n}$. The initial bound of Durfee,
 $\mu\ge n!p_g$, is for $\al=0$, $C_n=0$.
\li In section \ref{Sec.What.Is.Determined.By.Spectrum} we address the question "What is (not) determined by the spectrum?" In particular, the curve multiplicity is determined and in many cases the multiplicities of all the exceptional divisors of the minimal good resolution are determined.
 On the other hand, using recombination formulae it is immediate to construct lots of singularity sequences of distinct singularity types but with coinciding spectra, section  \ref{Sec.Sings.With.Coinciding.Spectra}.
\li In section \ref{Sec.Givental.Bound} we prove the Givental-type inequality for spectral values.
 Write $Spec(C,0)=a_0t^0+\suml_i a_i (t^{\al_i}+t^{-\al_i})$, where $0<\al_1<\al_2<\cdots<\al_k$, while $a_i\in\Z_{>0}$.
  We prove $\al_{r+i}+\al_{k-i}\le 1$, $i=0,1,\dots$, where $r\ge1$ is determined by $(C,0)$.

\

Finally, we remark that many \cp\ formulae hold in higher dimensions, for isolated hypersurface singularities. We hope to report on this in the subsequent paper.

\section{Preliminaries}

\subsection{The spectral pairs}
Let $(S,0)$ be a normal surface singularity which is a rational homology sphere. Consider a reduced curve singularity, $(C,0)\sset(S,0)$. Consider its embedded good resolution, see the diagram.
\\\parbox{14cm}{
Here
\li $U$ is a smooth surface, $\pi$ is proper and $U\smin\{E_i\}\isom(S,0)\smin\{0\}$.
\li $\{E_i\}$ are smooth projective curves, $|E_i\cap E_j|\le1$.
 $\tC\to(C,0)$ is the strict transform, $\tC$ is the multi-germ. All the intersections in $\tC\cup\{E_i\}$
 are normal crossings.
\li  $m_i=ord_{E_i}\pi^*(f)$
}\hspace{0.5cm}$\bM \tC&\sset&(U,\{E_i\})\\\downarrow&&\downarrow\pi\\(C,0)&\sset&(S,0)\\\downarrow&&\downarrow f\\\{0\}&\sset&(\C^1,0)\eM$

As $(S,0)$ is a rational homology sphere, the dual graph, $\Ga=\{E_i\}$, is a tree of $\P^1$'s.

Define the "chronology" on $\Ga$ as follows. Choose any irreducible component of $\{E_i\}$, call it the base point, $E_*$.
 For any other $E_i$ consider the shortest path from $E_i$ to $E_*$ in $\Ga$.
 (As $\Ga$ is a tree, this path is unique.)
 Denote the nearest neighbour of $E_i$ on this path by $E_{p(i)}$.
We often take as $E_*$ the first exceptional divisor that is born at the first blowup of $(S,0)$.

We recall the description of the spectral pairs,  $Spp(C,0)=\suml^\mu_{i=1}t^{[\xi_i,m_i]}$,
  from \cite{Schr.Stee.Stev.}, see also \cite{Thanh-Steenbrink-89} for the spectrum.

For each component $i\in\Ga\smin *$ introduce the following greatest common divisors:
\beq
\de_i=gcd(m_i,m_{p(i)}),\quad
r_i=\Bigg\{\ber 1:\ if\ \tC\cap E_i\neq\varnothing\\gcd(m_i,\{m_j\}_{E_j\cap E_i\neq\varnothing}):\ if\ \tC\cap E_i=\varnothing\eer.
\eeq

Using them we construct the following parts of the spectrum:
\beq\ber
a_i=\suml_{\substack{0<s<m_i\\m_i\nmid sr_i}}\Big(-1+|\tC\cap E_i|\frac{s}{m_i}+
\suml_{E_j\cap E_i\neq\varnothing}\{\frac{sm_j}{m_i}\}\Big)\Big(t^{[\frac{s}{m_i}-1,1]}+t^{[1-\frac{s}{m_i},1]}\Big)\\
b_i=\suml_{0<s<r_i}\Big(t^{[-\frac{s}{r_i},2]}+t^{[\frac{s}{r_i},0]}\Big),\quad
c_{i\neq *}=\suml_{0<s<\de_i}\Big(t^{[-\frac{s}{\de_i},2]}+t^{[\frac{s}{\de_i},0]}\Big),\quad
\eer\eeq
(Here $\{x\}$ is the fractional part of $x$.)
 Then, \cite[Theorem 2.1]{Schr.Stee.Stev.}:
\beq\label{Eq.Spectral.Pairs.Expression}
Spp(C,0)=\suml_{i\in\Ga}a_i+\suml_{*\neq i\in\Ga}(c_i-b_i)-b_*+(|\tC\cap E|-1)t^{[0,1]}
\eeq

\subsection{Some typical singularity types}
In this paper we work with embedded topological singularity types and associated invariants. Thus, by a curve singularity $(C,0)$, we mean any representative of the given singularity type.
 We often depict the singularity type by its (partial) resolution or by the dual graph.
\bee
\item Denote the ordinary multiple point (represented e.g. by $x^p=y^p$) of multiplicity $p$ by
\[\begin{picture}(60,10)(0,0)\linethickness{0.8mm}\curve(0,0,40,0)
\thinlines\freev{20}{-5}{10}\put(45,-3){$\tinyM pE_1$}
\put(10,-9){$\underbrace{\hspace{0.5cm}}_{p}$}
\end{picture}
\quad or,\ dually,\ by\quad
\begin{picture}(140,20)(-10,0)
\put(-2.5,1){$\bullet$}\claw{0}{2}{10}\put(-7,8){$\tinyM[pE_1]$}
\put(-9,-7){$\underbrace{\hspace{0.2cm}}_{p}$}
\end{picture}
\]
\

\item Denote the singularity type of $(x^p+y^p)(y^q-x^{2q})$ by
\[\begin{picture}(80,30)(0,0)\linethickness{0.8mm}\curve(10,0,50,0)
\thinlines\freev{20}{-5}{10}
\put(10,-9){$\underbrace{\hspace{0.5cm}}_{p}$}
\freeh{30}{20}{10}
\put(52,-2){$\tinyM(p+q)E_1$}
\linethickness{0.8mm}\curve(35,-10,35,30)
\put(45,18){$\Big\}\ q$}
\put(27,-20){$\tinyM (p+2q)E_2$}
\end{picture}
\quad\quad or,\ dually,\ by\quad
\begin{picture}(0,0)(-20,0)
\put(-2.5,1){$\bullet$}\put(-25,10){$\tinyM [(p+q)E_1]$}
\claw{0}{1}{8}\put(-9,-7){$\underbrace{\hspace{0.2cm}}_{p}$}
\curve(-1,3,25,3)
\put(25.5,1){$\bullet$}\put(20,10){$\tinyM [(p+2q)E_2]$}
\claw{28}{1}{8}\put(20,-7){$\underbrace{\hspace{0.3cm}}_{q}$}
\end{picture}
\]

\

\

\item The resolution of $(y^{p_1}+x^{p_1})(y^{p_2}-x^{2p_2})(y^{2p_3}-x^{3p_3})$ is given by
\[
\begin{picture}(130,60)(0,0)\linethickness{0.8mm}\curve(10,0,50,0)
\thinlines\freev{20}{-5}{10}
\put(10,-9){$\underbrace{\hspace{0.7cm}}_{p_1}$}
\put(52,-2){$\tinyM(p_1+p_2+2p_3)E_1$}
\linethickness{0.8mm}\curve(35,-10,35,45)
\thinlines\freeh{30}{20}{10}
\put(27,-20){$\tinyM (2p_1+3p_2+6p_3)E_3$}
\linethickness{0.8mm}\curve(30,40,75,40)
\thinlines\freev{55}{35}{10}
\put(45,45){$\overbrace{\hspace{0.7cm}}^{p_2}$}
\put(78,40){$\tinyM(2p_1+2p_2+3p_3)E_2$}
\end{picture}
\quad or\ dually\ \quad
\begin{picture}(40,0)(-20,0)
\curve(-1,3,52,3)
\put(-2.5,1){$\bullet$}\put(-20,10){$\tinyM [m_1E_1]$}
\claw{0}{1}{8}\put(-9,-7){$\underbrace{\hspace{0.2cm}}_{p_1}$}
\put(25.5,1){$\bullet$}\put(15,10){$\tinyM [m_3E_3]$}
\claw{28}{1}{8}\put(20,-7){$\underbrace{\hspace{0.3cm}}_{p_3}$}
\put(48,1){$\bullet$}\put(45,10){$\tinyM [m_2E_2]$}
\claw{50}{1}{8}\put(42,-7){$\underbrace{\hspace{0.3cm}}_{p_2}$}
\end{picture}
\hspace{2.5cm} \ber m_1=p_1+p_2+2p_3\\m_2=p_1+2p_2+3p_3\\m_3=2p_1+3p_2+6p_3\eer
\]

\

\item More generally, consider the singularity $\{(y^{p_1}+x^{p_1})\prodl^{k}_{i=2}(y^{(i-1)p_i}-x^{ip_i})=0\}\sset(\C^2,0)$. We call these singularities the  "intermediate types" and denote them by $C_{p_1,\dots,p_k}$.
The dual graph for the resolution of $C_{p_1,\dots,p_k}$ is:
\[\begin{picture}(120,20)(-20,0)
\curve(-1,3,123,3)

\put(-2.5,1){$\bullet$}\thinlines\claw{0}{1}{7} \put(-15,10){$\tinyM [m_1E_1]$}
\put(25,1){$\bullet$}\claw{27}{1}{7}\claw{27}{1}{7}\put(14,10){$\tinyM [m_kE_k]$}
\put(55,1){$\bullet$}\claw{57}{1}{7}\claw{27}{1}{7}\put(43,10){$\tinyM [m_{k-1}E_{k-1}]$}
\put(90,10){$\ldots$}
\put(120,1){$\bullet$}\claw{122}{1}{7}\claw{27}{1}{7}\put(115,10){$\tinyM [m_2E_2]$}

\put(-9,-7){$\underbrace{\hspace{0.2cm}}_{p_1}$}
\put(20,-7){$\underbrace{\hspace{0.3cm}}_{p_k}$}
\put(48,-7){$\underbrace{\hspace{0.3cm}}_{p_{k-1}}$}
\put(118,-7){$\underbrace{\hspace{0.3cm}}_{p_2}$}
\end{picture}
\hspace{2.5cm}
\ber m_1=p_1+p_2+2p_3+\cdots+(k-1)p_k\\m_2=m_1+\suml^k_{i=2}p_i\\\dots\\m_k=(k-1)m_1+p_2+2p_3+\cdots+(k-1)p_k=km_1-p_1\eer
\]
\eee
For further reference we record the Milnor numbers:
\bel\ \\
1. $\mu\Big((y^{p_1}+x^{p_1})\prodl^{k}_{i=2}(y^{(i-1)p_i}-x^{ip_i})\Big)=(p_1-1)^2+2p_1\suml^r_{i=2}(i-1)p_i+\suml_{i=2}^r(p_i(i-1)-1)(ip_i-1)+(1-r)+
2\suml_{2\le i<j}jp_ip_j$.
\\2. $\mu\Big((y^{p_1}+x^{p_1})\prodl^{k}_{i=2}(y^{(i-1)p_i}-x^{ip_i})\Big)-
\mu\Big(y^{p_1+\suml^r_{i=2}(i-1)p_i}-x^{p_1+\suml^r_{i=2}(i-1)p_i}\Big)
=1+\suml^r_{i=2}(p^2_i(i-1)-p_i)+2\suml^r_{2\le i<j}(i-1)p_ip_j$.
\eel
\bpr
We use the classical relation for $\de$-invariants, \cite[pg.206]{GLS}: $\de(\cupl_i(C_i,0))=\suml_i\de(C_i,0)+\suml_{1\le i<j}deg\Big((C_i,0),(C_j,0)\Big)$. Here the last term is the sum of the intersection numbers. To compute the Milnor number we use the relation $\de=\frac{\mu+r-1}{2}$, \cite{Buchweitz-Greuel-1980}.

In our case: $(C_1,0)=\{x^{p_1}=y^{p_1}\}$, while $(C_{i>1},0)=\{y^{(i-1)p_i}=x^{ip_i}\}$.
 Therefore:
\beq
\mu(x^{p_1}-y^{p_1})=(p_1-1)^2,\quad \mu(y^{(i-1)p_i}-x^{ip_i})=\Big((i-1)p_i-1\Big)(ip_i-1)\quad and\quad
 deg\Big((C_1,0),(C_i,0)\Big)=p_1(i-1)p_i.
\eeq
 It remains to compute
  $deg\Big((C_i,0),(C_j,0)\Big)=p_ip_jdeg\Big(y^{(i-1)}-x^{i},y^{(j-1)}-x^{j}\Big)$, for $1<i<j$.
This intersection multiplicity can be computed by restricting the ideal $(y^{(j-1)}-x^{j})$ to the local ring of $y^{(i-1)}=x^{i}$.
 Namely, inside the ring $\k[[t^{i-1},t^i]]$ one considers the ideal generated by $t^{ij-j}-t^{ij-i}$.
  We want to compute $dim_\k\quotient{\k[[t^{i-1},t^i]]}{(t^{ij-j}-t^{ij-i})}$.

The conductor of the local ring is $t^{(i-1)(i-2)}$ and the basis (as a $\k$-vector space) is: $1,t^{i-1},t^i,t^{2i-1},\dots,\{t^{(i-1)(i-2)+k}\}_{k=0,1\dots}$. Therefore the conductor of the ideal is
 $t^{(i-1)(i-2)+ij-j}$, while the $\k$-basis of the ideal is
\beq
t^{ij-j}-t^{ij-i},t^{i-1}(t^{ij-j}-t^{ij-i}),t^i(t^{ij-j}-t^{ij-i}),\dots,\{t^{(i-1)(i-2)+ij-j+k}\}_{k=0,1,\dots}.
\eeq
Note that the gaps of the ideal are in the natural bijection with the gaps of the ring.
 Finally,
\beq
dim_\k\quotient{\k[[t^{i-1},t^i]]}{(t^{ij-j}-t^{ij-i})}=ij-j-\sharp(gaps\ of\ \k[[t^{i-1},t^i]])+
\sharp(gaps\ of\ the\ ideal)=(i-1)j
\eeq
Thus for $i<j$: $deg\Big((C_i,0),(C_j,0)\Big)=p_ip_j(i-1)j$. Substitute all this into the relation for $\de$-invariant, to get:
\beq
\frac{\mu+\suml^r_{i=1}p_i-1}{2}=
\frac{(p_1-1)^2+\suml^r_{i=2}((i-1)p_i-1)(ip_i-1)+\suml^r_{i=1}(p_i-1)}{2}+p_1\suml^r_{i=2}(i-1)p_i+
\suml^r_{1<i<j}p_ip_j(i-1)j,
\eeq
proving the first statement.

The second statement follows by subtraction, $\mu-(p_1+\sum(i-1)p_i-1)^2$.
\epr

\section{Additivity of the spectrum}
\bel\label{Thm.Spectrum.Additive.Formula}
Consider a (good) embedded resolution $\tC\sset(U,\{E_i\})\to(S,0)\supset(C,0)$. Then
\[Spec(C,0)=\suml_{i\in\Ga}\suml_{|k|<m_i}d_{|k|,i}t^{\frac{k}{m_i}},\quad where\
 d_{k,i}=-1+|\tC\cap E_i|\frac{m_i-k}{m_i}+\suml_{\substack{E_j\cap E_i\neq\varnothing\\j\neq p(i)}}\{-\frac{km_j}{m_i}\}+\underbrace{1-\{\frac{km_{p(i)}}{m_i}\}}_{only\ if\ i\neq*}.
\]
\eel
\bpr
The proof is just a rearrangement of the original expression in equation (\ref{Eq.Spectral.Pairs.Expression}).
We need only the spectral values, thus, ignoring the second arguments in $t^{[\xi,q]}$, we get $Spec(C,0)=\suml_{i\in\Ga}Spec_{E_i}$, where each $Spec_{E_i}$ consists of the contributions from $E_i$.

First note that $b_i\neq0$ only if $r_i>1$ i.e. $\tC\cap E_i=\varnothing$. But in this case
 $m_i\nmid sr_i$ iff $\frac{s}{m_i}=\frac{s'}{r_i}$ for some $s'\in\N$ iff $\{\frac{sm_j}{m_i}\}=0$.
Therefore:
\beq
\suml_{\substack{0<s<m_i\\m_i\nmid sr_i}}\Big(\dots\Big)-b_i=
\suml_{0<s<m_i}\Big(\dots\Big).
\eeq
So, we get: $Spec_{E_i}=\suml_{0<s<m_i}\Big(\dots\Big)+c_i+|\tC\cap E_i|t^0-\underbrace{t^0}_{if\ i=*}$, with the convention
$c_i=\Bigg\{\ber 0:\ i=*\\\suml_{0<s<\de_i}\Big(t^{-\frac{s}{\de_i}}+t^{\frac{s}{\de_i}}\Big)\eer$.

We claim that this expression coincides with the stated one.
In the case $i=*$ the coincidence is directly seen.
So we consider the case $i\neq*$. By renaming $s\to m_i-s$ we get:
\beq
Spec_{E_i}=\suml_{0<s<m_i}\Big(-1+|\tC\cap E_i|\frac{m_i-s}{m_i}+\suml_{E_j\cap E_i\neq\varnothing}\{-\frac{sm_j}{m_i}\}\Big)\Big(t^{\frac{s}{m_i}}+t^{-\frac{s}{m_i}}\Big)+
\suml_{0<s<\de_i}(t^{\frac{s}{\de_i}}+t^{-\frac{s}{\de_i}})+|\tC\cap E_i|t^0
\eeq
Note that $\{-\frac{sm_{p(i)}}{m_i}\}=0$ iff $\frac{s}{m_i}=\frac{s'}{\de_i}$, for some (unique) $\de_i> s'\in\N$.
 Further, note that the function $h(x)=\Big\{\ber \{-x\}:\ x\not\in\Z\\1:x\in\Z\eer$ can be written in the form $h(x)=1-\{x\}$. Therefore
\beq
Spec_{E_i}=\suml_{0<s<m_i}\Big(-1+|\tC\cap E_i|\frac{m_i-s}{m_i}+
\suml_{\substack{E_j\cap E_i\neq\varnothing\\j\neq p(i)}}\{-\frac{sm_j}{m_i}\}+\underbrace{1-\{\frac{sm_{p(i)}}{m_i}\}}_{only\ if\ i\neq*}\Big)
(t^{\frac{s}{m_i}}+t^{\frac{-s}{m_i}})+|\tC\cap E_i|t^0.
\eeq
Now the statement follows.
\epr

\bex\label{Ex.OMP.Basic.Types.Spectrum}
1. For the \omp, $(C,0)=\{x^m=y^m\}$, we get: $Spec(C,0)=\suml_{|k|<m}(m-|k|-1)t^{\frac{k}{m}}$.
 The maximal spectral value is  $\al_{max}=\frac{m-2}{m}$.

2. Another particularly important case is $C_{p,q}=\{(x^p+y^p)(y^q-x^{2q})=0\}$, $q\ge2$. (For $q=0$ or $q=1$ this
 is just an \omp.) We call this singularity "of the basic type". Lemma \ref{Thm.Spectrum.Additive.Formula} gives:
\beq
Spec(C_{p,q})=\suml_{|k|<p+2q}(q-\lceil\frac{q|k|}{p+2q}\rceil)t^{\frac{k}{p+2q}}+
\suml_{|k|<p+q}(p-1-\lfloor\frac{p|k|}{p+q}\rfloor)t^{\frac{k}{p+q}}
\eeq
We need the corresponding maximal spectral value, $\al_{max}$. Present the spectrum in the form
\beq
\suml_{1\le k<p+2q}\lfloor\frac{qk}{p+2q}\rfloor (t^{1-\frac{k}{p+2q}}+t^{-1+\frac{k}{p+2q}})+
\suml_{1\le k<p+q}(\lceil\frac{pk}{p+q}\rceil-1)(t^{1-\frac{k}{p+q}}+t^{-1+\frac{k}{p+q}})+t^0(q+p-1).
\eeq
Then we get: if $p>q$ then $\al_{max}=1-\frac{2}{p+q}$, if $p\le q$ then $\al_{max}=1-\frac{3}{p+2q}$.
  (To check the case $p\le q$ one distinguishes between $p=0$ and $p>0$.)
\eex

\section{Cut-and-paste formulae}
\ \\
\parbox{13cm}{\subsection{}
Suppose after the i'th blowup one gets a smooth surface with the configuration as on the first picture.
Here $D$ denotes some collection of curve singularities, while the blackbox denotes the
 rest of exceptional divisors and the remaining part of the strict transform, $\tC\smin D$.

Consider the associated singularity types, defined by the partial resolutions.
On the  second  picture the curve singularities of $D$ are replaced by several smooth curve germs
 that intersect $E_i$ normally.
  On the third picture the rest of exceptional divisors is replaced by a collection of smooth germs.
}
\begin{picture}(0,0)(-50,-30)
\curve(0,0,40,0)\put(20,0){$\curlyvee$}\put(30,0){$\curlyvee$}
\put(20,8){$\overbrace{}$}  \put(25,17){$\tinyM D$}
\put(-5,-3){\blackboxh}\put(-28,-2){$\tinyM m_i E_i$}
\end{picture}
\begin{picture}(10,0)(-50,0)
\curve(0,0,40,0)\freev{30}{-5}{10}
\put(20,8){$\overbrace{}$}  \put(20,17){$\tinyM deg(D\cap E_i)$}
\put(-5,-3){\blackboxh}\put(-28,-2){$\tinyM m_i E_i$}
\end{picture}
\begin{picture}(0,0)(-30,40)
\curve(0,0,40,0)\freev{10}{-5}{10}\freev{30}{-5}{10}
\put(5,8){$\overbrace{\quad\quad\quad}$}  \put(15,17){$\tinyM m_i$}
\put(-25,-2){$\tinyM m_i E_i$}
\end{picture}
\\\\

\bthe\label{Thm.Spectrum.Cp.Formulas}
1.
\begin{picture}(0,0)(-30,0)
\curve(0,0,40,0)\put(20,0){$\curlyvee$}\put(30,0){$\curlyvee$}
\put(20,8){$\overbrace{}$}  \put(25,17){$\tinyM D$}
\put(-5,-3){\blackboxh}\put(-28,-2){$\tinyM m_i E_i$}
\put(50,0){$-$}
\end{picture}
\begin{picture}(0,0)(-120,0)
\curve(0,0,40,0)\freev{30}{-5}{10}
\put(20,8){$\overbrace{}$}  \put(20,17){$\tinyM deg(D\cap E_i)$}
\put(-5,-3){\blackboxh}\put(-28,-2){$\tinyM m_i E_i$}
\put(50,0){$\equiv$}
\end{picture}
\begin{picture}(0,0)(-210,0)
\curve(0,0,40,0)\put(20,0){$\curlyvee$}\put(30,0){$\curlyvee$}
\put(20,8){$\overbrace{}$}  \put(25,17){$\tinyM D$}
\freev{10}{-5}{10}\put(-25,-2){$\tinyM m_i E_i$}
\put(0,-7){$\underbrace{}$}\put(-7,-20){$\tinyM m_i-deg(D\cap E_i)$}
\put(50,0){$-$}
\end{picture}
\begin{picture}(0,0)(-300,0)
\curve(0,0,40,0)\freev{10}{-5}{10}\freev{30}{-5}{10}
\put(5,8){$\overbrace{\quad\quad\quad}$}  \put(15,17){$\tinyM m_i$}
\put(-25,-2){$\tinyM m_i E_i$}
\end{picture}\hspace{8cm}
\\\\\\

or dually:
\begin{picture}(0,0)(-10,-3)
\curve(0,0,20,0)\thinlines
\put(10,-2){$\blackboxh$}
\put(-2.5,-2.5){$\bullet$}
\put(-4,-7){$\curlywedge$}\put(-4,-15){$\tinyM D$}
\put(-5,8){$\tinyM [m_iE_i]$}\put(35,-3){$=$}
\end{picture}
\begin{picture}(0,0)(-60,-3)
\curve(0,0,20,0)
\put(10,-2){$\blackboxh$}
\put(-2.5,-2.5){$\bullet$}\claw{0}{1}{10}
\put(-18,-18){${\tinyM\tinyT deg(D\cap E_i)}$}
\put(-10,8){$\tinyM [m_iE_i]$}\put(35,-3){$+$}
\end{picture}
\begin{picture}(0,0)(-110,-3)
\put(-2.5,-2.5){$\bullet$}\claw{0}{1}{10}\put(-4,2){$\curlyvee$}
\put(-18,-18){${\tinyM\tinyT deg(D\cap E_i)}$}
\put(-10,10){$\tinyM [m_iE_i]$}\put(12,-3){$-$}
\end{picture}
\begin{picture}(0,0)(-150,-3)
\put(-2.5,-2.5){$\bullet$}
\claw{0}{1}{10}
\curve(0,0,0,10)\curve(0,0,-5,10)\curve(0,0,5,10)
\put(-18,-18){${\tinyM\tinyT deg(D\cap E_i)}$}
\put(-10,15){$\tinyM [m_iE_i]$}
\end{picture}
\vspace{1cm}

2. For any reduced curve singularity the spectrum is a linear combination of the "intermediate" types:

\

$(C,0)\equiv\suml_{\{m_1,\dots,m_k\}}$
\begin{picture}(0,0)(-10,-3)
\curve(0,0,135,0)
\put(-2.5,-2.5){$\bullet$}\put(37.5,-2.5){$\bullet$}\put(77.5,-2.5){$\bullet$}\put(132.5,-2.5){$\bullet$}
\claw{0}{1}{10}\claw{40}{1}{10}\claw{80}{1}{10}\claw{135}{1}{10}
\put(-10,5){$\tinyM[m_1E_1]$}\put(30,5){$\tinyM[m_kE_k]$}\put(65,5){$\tinyM[m_{k-1}E_{k-1}]$}
\put(125,5){$\tinyM[m_2E_2]$}
\end{picture}
\begin{picture}(0,0)(-165,-3)
\put(-5,-3){$-\suml_{m_j}$}\put(27.5,-2.5){$\bullet$}\claw{30}{1}{10}
\put(20,5){$\tinyM[m_jE_j]$}
\end{picture}
\

Here the first sum consists of the singularities of types $(y^{p_1}-x^{p_1})\prod^{k}_{i=2}(y^{(i-1)p_i}-x^{ip_i})$, while the second sum consists of the ordinary multiple points.
The number of summands in the first sum is by one bigger than that in the second sum.

3. If all the branches of  $(C,0)$ are smooth then the decomposition simplifies:

\

$(C,0)\equiv\suml_{\{m_1,m_2\}}$
\begin{picture}(0,0)(-10,-3)
\curve(0,0,40,0)
\put(-2.5,-2.5){$\bullet$}\put(37.5,-2.5){$\bullet$}
\claw{0}{1}{10}\claw{40}{1}{10}
\put(-10,5){$\tinyM[m_1E_1]$}\put(30,5){$\tinyM[m_2E_2]$}
\end{picture}
\begin{picture}(0,0)(-75,-3)
\put(-5,-3){$-\suml_{m_i}$}\put(27.5,-2.5){$\bullet$}\claw{30}{1}{10}
\put(20,5){$\tinyM[m_jE_j]$}
\end{picture}
\vspace{0.5cm}

4. Finally, the intermediate types decompose as follows:
\beqm\label{Eq.Spec.Decomposition.of.Intermediate.Types}
\begin{picture}(0,0)(-10,-3)
\curve(0,0,125,0)
\put(-2.5,-2.5){$\bullet$}\put(37.5,-2.5){$\bullet$}\put(77.5,-2.5){$\bullet$}\put(122.5,-2.5){$\bullet$}
\claw{0}{1}{10}\claw{40}{1}{10}\claw{80}{1}{10}\claw{125}{1}{10}
\put(-2,-15){$\tinyM p_1$}\put(38,-15){$\tinyM p_k$}\put(75,-15){$\tinyM  p_{k-1}$}\put(100,-10){$\dots$}
\put(123,-15){$\tinyM p_2$}
\put(-10,5){$\tinyM[m_1E_1]$}\put(30,5){$\tinyM[m_kE_k]$}\put(65,5){$\tinyM[m_{k-1}E_{k-1}]$}
\put(115,5){$\tinyM[m_2E_2]$}\put(140,-3){$\equiv$}
\end{picture}
\begin{picture}(0,0)(-170,-3)
\put(-2.5,-2.5){$\bullet$}\put(-10,5){$\tinyM [m_2E_1]$}\put(15,-3){$+$}
\end{picture}
\begin{picture}(0,0)(-240,-3)
\put(-40,-3){$\suml^{k-1}_{i=2}\!\Big($}
\curve(0,0,40,0)\put(-2.5,-2.5){$\bullet$}\put(37.5,-2.5){$\bullet$}
\claw{0}{1}{10}\claw{40}{1}{10}
\put(-24,-20){$\tinyM m_i-m_1-\suml_{j>i}p_i$}\put(30,-17){$\tinyM m_1+\suml_{j>i}p_i$}
\put(-20,7){$\tinyM [m_iE_1]$}\put(25,7){$\tinyM [m_{i+1}E_2]$}\put(60,-3){$-$}
\end{picture}
\begin{picture}(0,0)(-330,-3)
\put(-2.5,-2.5){$\bullet$}\claw{0}{1}{10}\put(-20,7){$\tinyM [m_{i+1}E_1]$}
\put(-5,-17){$\tinyM m_{i+1}$}
\put(15,-3){$\Big)$}\put(22,-3){$+$}
\end{picture}
\\
\\
\\
\begin{picture}(0,0)(400,-3)
\put(-20,0){$+$}\curve(0,0,40,0)
\put(-2.5,-2.5){$\bullet$}\put(37.5,-2.5){$\bullet$}
\claw{0}{1}{10}\claw{40}{1}{10}
\put(-2,-15){$\tinyM p_1$}\put(25,-15){$\tinyM m_1-p_1$}
\put(-10,5){$\tinyM [m_1E_1]$}\put(20,5){$\tinyM [(2m_1-p_1)E_2]$}
\put(65,-3){$-$}
\end{picture}
\begin{picture}(0,0)(310,-3)
\put(-2.5,-2.5){$\bullet$}\claw{0}{1}{10}\put(-15,-17){$\tinyM 2m_1-p_1$}
\put(15,-3){$+$}
\end{picture}
\begin{picture}(0,0)(230,-3)
\put(-35,-3){$\suml^{k-1}_{i=2}\!\Big($}
\put(-2.5,-2.5){$\bullet$}\claw{0}{1}{10}
\put(-12,7){$\tinyM [(im_1-p_1)E_1]$}
\put(25,-3){$-$}
\end{picture}
\begin{picture}(0,0)(160,-3)
\curve(0,0,70,0)\put(-2.5,-2.5){$\bullet$}\put(67.5,-2.5){$\bullet$}
\claw{0}{1}{10}\claw{70}{1}{10}
\put(-14,-17){$\tinyM (i-1)m_1-p_1$}\put(65,-17){$\tinyM m_1$}
\put(-20,7){$\tinyM [(im_1-p_1)E_1]$}\put(35,7){$\tinyM [(i+1)m_1-p_1)E_2]$}
\put(100,-3){$\Big)$}
\end{picture}
\end{multline}
\ethe
Thus the spectrum of any plane curve singularity decomposes in the a linear combination of the spectra of basic types, $C_{p,q}=\{(x^p+y^p)(y^q-x^{2q})=0\}\sset(\C^2,0)$.
\\\bpr
1. Use the presentation of spectrum obtained in lemma \ref{Thm.Spectrum.Additive.Formula}.
 As the expression is additive, $Spec(C,0)=\suml_i Spec_{E_i}$, it is enough to check in the formula each $Spec_{E_j}$.
 As the \cp\ operation occurs only at $E_i$, we should check only $Spec_{E_i}$.

Finally, the expression for $Spec_{E_i}$ is additive in $\tC\cap E_i$ and $E_j\cap E_i$.

2. The proof is an easy induction on the minimal number of blowups needed to resolve the singularity.
 Blowup $(C,0)\sset(\C^2,0)$ once, let $E_1$ be the exceptional divisor, and $\tC^1$ the strict transform.
 If the intersection $E_1\cap\tC^1$ is non-transverse in at least two (distinct) points then apply equation (\ref{Eq.Spectral.Relation.Simplest.Case}).
 (It is a particular case of relations of part 1.) Therefore we can assume that either $\tC^1$, $E_1$ intersect  transversely everywhere, i.e. $(C,0)$ was an ordinary multiple point, or the intersection  is non-transverse precisely at one point. In the later case blowup at this point. Let $E_2$ be the new component, let $\tC^2$ be the current strict transform of $(C,0)$, denote the strict transform of $E_1$ by $\tE_1$.

By the construction, the intersection of $\tC^2$ with $\tE_1\smin E_2$ is transverse.
 If the intersection of $\tC^2$ with $\tE_1\cup E_2$ is transverse then we get the statement.
If the intersection of $\tC^2$ with $E_2\smin\tE_1$ is non-transverse, apply part 1. Then one get a singularity whose resolution length is smaller than that of $(C,0)$, hence one uses the induction assumption.
 The remaining case is: the only non-transversality occurs at the point $\tE_1\cap E_2$. In this case blowup this point
  $\tE_1\cap E_2$ and continue in the same way.

3. If all the branches of $(C,0)$ are smooth, then no branch of the strict transform passes through the intersection point $E_i\cap E_j$. Thus we do not need to blowup at the point $E_i\cap E_j$.

4.  {\bf Step 1.} We prove the elementary decomposition:
\beqm\label{Eq.Spec.Elem.Decomp.Intermediate.Types}
\begin{picture}(0,0)(-10,-3)
\curve(0,0,135,0)
\put(-2.5,-2.5){$\bullet$}\put(37.5,-2.5){$\bullet$}\put(77.5,-2.5){$\bullet$}\put(132.5,-2.5){$\bullet$}
\claw{0}{1}{10}\claw{40}{1}{10}\claw{80}{1}{10}\claw{135}{1}{10}
\put(-2,-15){$\tinyM p_1$}\put(38,-15){$\tinyM p_k$}\put(75,-15){$\tinyM p_{k-1}$}\put(133,-15){$\tinyM p_2$}
\put(-10,5){$\tinyM[m_1E_1]$}\put(30,5){$\tinyM[m_kE_k]$}\put(65,5){$\tinyM[m_{k-1}E_{k-1}]$}
\put(125,5){$\tinyM[m_2E_2]$}\put(160,-3){$-$}
\end{picture}
\begin{picture}(0,0)(-200,-3)
\curve(0,0,105,0)
\put(-2.5,-2.5){$\bullet$}\put(37.5,-2.5){$\bullet$}\put(102.5,-2.5){$\bullet$}
\claw{0}{1}{10}\claw{40}{1}{10}\claw{105}{1}{10}
\put(-10,-18){$\tinyM p_k+m_1$}\put(38,-18){$\tinyM p_{k-1}$}\put(70,-10){$\dots$}
\put(102,-18){$\tinyM  p_2$}
\put(-15,5){$\tinyM[m_kE_k]$}
\put(70,5){$\dots$}
\put(20,5){$\tinyM [m_{k-1}E_{k-1}]$}
\put(90,5){$\tinyM[m_2E_2]$}\put(125,-3){$\equiv$}
\end{picture}
\\\\\\\hspace{3cm}
\begin{picture}(0,0)(450,-3)
\put(-20,-3){$\equiv$}\curve(0,0,105,0)
\put(-2.5,-2.5){$\bullet$}\put(37.5,-2.5){$\bullet$}\put(102.5,-2.5){$\bullet$}
\claw{0}{1}{10}\claw{40}{1}{10}\claw{105}{1}{10}
\put(-2,-15){$\tinyM p_1$}\put(35,-18){$\tinyM\tp_{k-1}$}\put(65,-10){$\dots$}\put(103,-18){$\tinyM \tp_2$}
\put(-10,5){$\tinyM [m_1E_1]$}\put(20,5){$\tinyM [(m_k-m_1)E_{k-1}]$}\put(90,5){$\tinyM [\tilde{m}_2 E_2]$}
\put(120,-3){$-$}
\end{picture}
\begin{picture}(0,0)(300,-3)
\curve(0,0,125,0)
\put(-2.5,-2.5){$\bullet$}\put(52.5,-2.5){$\bullet$}\put(122.5,-2.5){$\bullet$}
\claw{0}{1}{10}\claw{55}{1}{10}\claw{125}{1}{10}
\put(-10,-17){$\tinyM \tp_{k-1}+m_1$}\put(48,-17){$\tinyM\tp_{k-2}$}\put(95,-10){$\dots$}\put(95,5){$\dots$}
\put(123,-17){$\tinyM \tp_2$}
\put(-20,5){$\tinyM[(m_k-m_1)E_{k-1}]$}\put(45,5){$\tinyM[m_{k-2}E_{k-2}]$}\put(115,5){$\tinyM [\tilde{m}_2E_2]$}
\put(140,-3){$+$}
\end{picture}
\begin{picture}(0,0)(135,-3)
\put(-2.5,-2.5){$\bullet$}\claw{0}{1}{10}\put(-20,-17){$\tinyM m_k-m_1$}
\put(15,-3){$-$}
\end{picture}
\begin{picture}(0,0)(95,-3)
\curve(0,0,80,0)\put(-2.5,-2.5){$\bullet$}\put(77.5,-2.5){$\bullet$}
\claw{0}{1}{10}\claw{80}{1}{10}
\put(-14,-17){$\tinyM m_k-2m_1$}\put(75,-17){$\tinyM m_1$}
\put(-20,7){$\tinyM [(m_k-m_1)E_1]$}\put(60,7){$\tinyM [m_kE_2]$}
\end{picture}
\end{multline}
\\\\Here $m_1=p_1+\suml^k_{i=2}(i-1)p_i$, while $\{\tp_i\}_{i=1,\dots,k-1}$ are some additional multiplicities  satisfying $m_1=p_1+\suml^{k-1}_{i=2}(i-1)\tp_i$,
 and such that the corresponding graphs are the (dual) resolution graphs of some reduced plane curve singularities.

In terms of the representatives of singularity types this equality can be written as:
\beqm\label{Eq.Spec.Elem.Decomp.Formula}
Spec\Bigg((y^{p_1}-x^{p_1})\prodl^{k}_{i=2}(y^{(i-1)p_i}-x^{ip_i})\Bigg)-
Spec\Bigg((y^{p_k+m}-x^{(k-1)(p_k+m)})\prodl^{k-1}_{i=2}(y^{p_i}-x^{(i-1)p_i})\Bigg)=\\
=Spec\Bigg((y^{p_1}-x^{p_1})\prodl^{k-1}_{i=2}(y^{(i-1)\tp_i}-x^{i\tp_i})\Bigg)
-Spec\Bigg((y^{\tp_{k-1}+m}-x^{(k-2)(\tp_{k-1}+m)})\prodl^{k-2}_{i=2}(y^{\tp_i}-x^{(i-1)\tp_i})\Bigg)+\\
+Spec\Bigg(x^{(k-1)m-p_1}-y^{(k-1)m-p_1}\Bigg)-Spec\Bigg((x^{m}-y^{m})(y^{(k-2)m-p_1}-x^{2(k-2)m-2p_1})\Bigg)
\end{multline}

\parbox{12cm}{
The resolution of the singularity $\Big\{(y^{p_1}-x^{p_1})\prodl^{k}_{i=2}(y^{(i-1)p_i}-x^{ip_i})=0\Big\}$ is on the right.
 Here $m_1=p_1+\suml^k_{i=2}(i-1)p_i$, $m_2=m_1+\suml^k_{i=2}p_i$, $m_3=2m_1+p_2+2\suml^k_{i=2}p_i$, \dots
 $m_k=(k-1)m_1+p_2+2p_3+\cdots+(k-1)p_k=km-p_1$.

We should compare the sums of spectra in both parts of equation (\ref{Eq.Spec.Elem.Decomp.Formula}).
}
\begin{picture}(0,0)(-30,0)
\linethickness{0.8mm}\curve(0,0,30,0)\thinlines\freev{10}{-5}{10}
\put(-22,-2){$\tinyM mE_1$}\put(8,-10){$\tinyM p_1$}
\linethickness{0.8mm}\curve(25,-10,25,30)\thinlines\freeh{20}{15}{10}
\put(20,-18){$\tinyM m_kE_k$}\put(33,10){$\tinyM p_k$}
\linethickness{0.8mm}\curve(20,28,50,28)\thinlines\freev{35}{23}{10}
\put(-23,30){$\tinyM m_{k-1}E_{k-1}$}\put(30,38){$\tinyM p_{k-1}$}
\put(55,28){$\dots\dots$}
\linethickness{0.8mm}\curve(100,20,100,50)\thinlines\freeh{95}{40}{10}
\put(95,12){$\tinyM m_2E_2$}\put(110,40){$\tinyM p_2$}
\end{picture}

Use the presentation of lemma \ref{Thm.Spectrum.Additive.Formula}, $Spec(C,0)=\sum_{i\in\Ga}Spec_{E_i}$.
Note that there is the natural correspondence between $\{E_i\}$ of all the participants.
 (It is helpful to draw the actual resolutions.)
 For each $E_i$ of each participant the contribution consists of the part $|\tC\cap E_i|\frac{m_i-k}{m_i}t^{\frac{k}{m_i}}$ and the part related to all the vertices, $\{E_i\cap E_j\}$. By direct check: the parts coming from
$|\tC\cap E_i|\frac{m_i-k}{m_i}t^{\frac{k}{m_i}}$ coincide.

It remains to check the parts coming from the vertices. Again, the comparison of the resolution shows that most of the vertices cancel, we have to compare only the following contributions:\\
\beq
\begin{picture}(0,0)(130,0)
\linethickness{0.8mm}\curve(0,0,20,0)\curve(15,-5,15,15)
\put(-22,-2){$\tinyM m_1E_1$}\put(8,-13){$\tinyM m_kE_k$}\put(25,-3){$-$}
\curve(40,-5,40,15)\put(35,-13){$\tinyM m_k E_{k-1}$}\put(60,0){$\overset{?}{\equiv}$}
\end{picture}
\begin{picture}(0,0)(20,0)
\linethickness{0.8mm}\curve(0,0,20,0)\curve(15,-5,15,15)
\put(-22,-2){$\tinyM m_1E_1$}\put(-10,-13){$\tinyM (m_k-m_1)E_{k-1}$}\put(35,-3){$-$}
\curve(60,-5,60,15)\put(45,-13){$\tinyM (m_k-m_1)E_{k-2}$}\put(80,-3){$+$}
\end{picture}
\begin{picture}(0,0)(-120,0)\linethickness{0.8mm}
\curve(0,-5,0,15)\put(-22,-13){$\tinyM (m_k-m_1)E_1$}
\put(20,-3){$-$}
\curve(35,0,55,0)\curve(40,-5,40,15)
\put(30,-13){$\tinyM m_kE_2$}
\put(60,0){$\tinyM (m_k-m_1)E_1$}
\end{picture}
\eeq
\\

(Note: the strict transform $\tC$ and the rest of exceptional divisors are not drawn here.)

We write the formulae corresponding to the pictures. In all the cases we choose $*=2$.
\beqm
\Bigg(\suml_{|s|<m_k}t^{\frac{s}{m_k}}(-1+1+\{-\frac{s m}{m_k}\})+
\suml_{|s|<m}t^{\frac{s}{m}}(-1+1-\{\frac{s m_k}{m}\})\Bigg)
-\Bigg(\suml_{|s|<m_k}t^{\frac{s}{m_k}}(-1+1)\Bigg)
\stackrel{?}{=}\\\stackrel{?}{=}
\Bigg(\suml_{|s|<m_k-m}t^{\frac{s}{m_k-m}}(-1+1+\{-\frac{s m}{m_k-m}\})+
\suml_{|s|<m}t^{\frac{s}{m}}(-1+1-\{\frac{s(m_k-m)}{m}\})\Bigg)-\Bigg(\suml_{|s|<m_k-m}t^{\frac{s}{m_k-m}}(-1+1)\Bigg)
+\\+
\Bigg(\suml_{|s|<m_k-m}t^{\frac{s}{m_k-m}}(-1)\Bigg)
-\Bigg(\suml_{|s|<m_k}t^{\frac{s}{m_k}}(-1+1-\{\frac{s(m_k-m)}{m_k}\})
+\suml_{|s|<m_k-m}t^{\frac{s}{m_k-m}}(-1+\{-\frac{s m_k}{m_k-m}\})\Bigg)
\end{multline}

By direct check: everything cancels.

{\bf Step 2.} Now we use formulae (\ref{Eq.Spec.Elem.Decomp.Intermediate.Types}), (\ref{Eq.Spec.Elem.Decomp.Formula}) recursively:
\beqm
\begin{picture}(0,0)(-10,-3)
\curve(0,0,135,0)
\put(-2.5,-2.5){$\bullet$}\put(37.5,-2.5){$\bullet$}\put(77.5,-2.5){$\bullet$}\put(132.5,-2.5){$\bullet$}
\claw{0}{1}{10}\claw{40}{1}{10}\claw{80}{1}{10}\claw{135}{1}{10}
\put(-2,-15){$\tinyM p_1$}\put(38,-15){$\tinyM p_k$}\put(75,-15){$\tinyM p_{k-1}$}\put(133,-15){$\tinyM p_2$}
\put(-10,5){$\tinyM[m_1E_1]$}\put(30,5){$\tinyM[m_kE_k]$}\put(65,5){$\tinyM[m_{k-1}E_{k-1}]$}
\put(125,5){$\tinyM[m_2E_2]$}\put(160,-3){$-$}
\end{picture}
\begin{picture}(0,0)(-200,-3)
\curve(0,0,105,0)
\put(-2.5,-2.5){$\bullet$}\put(37.5,-2.5){$\bullet$}\put(102.5,-2.5){$\bullet$}
\claw{0}{1}{10}\claw{40}{1}{10}\claw{105}{1}{10}
\put(-10,-18){$\tinyM p_k+m_1$}\put(38,-18){$\tinyM p_{k-1}$}\put(70,-10){$\dots$}
\put(102,-18){$\tinyM  p_2$}
\put(-15,5){$\tinyM[m_kE_k]$}
\put(70,5){$\dots$}
\put(20,5){$\tinyM [m_{k-1}E_{k-1}]$}
\put(90,5){$\tinyM[m_2E_2]$}\put(125,-3){$\equiv$}
\end{picture}
\\\\\\\hspace{3cm}
\begin{picture}(0,0)(450,-3)
\put(-20,-3){$\equiv$}\curve(0,0,40,0)
\put(-2.5,-2.5){$\bullet$}\put(37.5,-2.5){$\bullet$}
\claw{0}{1}{10}\claw{40}{1}{10}
\put(-2,-15){$\tinyM p_1$}\put(25,-15){$\tinyM m_1-p_1$}
\put(-10,5){$\tinyM [m_1E_1]$}\put(20,5){$\tinyM [(m_k-(k-2)m_1)E_2]$}
\put(65,-3){$-$}
\end{picture}
\begin{picture}(0,0)(350,-3)
\put(-2.5,-2.5){$\bullet$}\claw{0}{1}{10}\put(-15,-17){$\tinyM m_k-(k-2)m_1$}
\put(15,-3){$+$}
\end{picture}
\begin{picture}(0,0)(260,-3)
\put(-50,-3){$\suml^{k-2}_{j=1}\!\Big($}
\put(-2.5,-2.5){$\bullet$}\claw{0}{1}{10}\put(-25,7){$\tinyM [(m_k-jm_1)E_1]$}
\put(25,-3){$-$}
\end{picture}
\begin{picture}(0,0)(205,-3)
\curve(0,0,70,0)\put(-2.5,-2.5){$\bullet$}\put(67.5,-2.5){$\bullet$}
\claw{0}{1}{10}\claw{70}{1}{10}
\put(-14,-17){$\tinyM m_k-(j-1)m_1$}\put(65,-17){$\tinyM m_1$}
\put(-20,7){$\tinyM [(m_k-jm_1)E_1]$}\put(35,7){$\tinyM [m_k-(j-1)m_1)E_2]$}
\put(100,-3){$\Big)$}
\end{picture}
\end{multline}

\

 Finally, present
\beq
\begin{picture}(0,0)(220,-3)
\curve(0,0,105,0)
\put(-2.5,-2.5){$\bullet$}\put(37.5,-2.5){$\bullet$}\put(102.5,-2.5){$\bullet$}
\claw{0}{1}{10}\claw{40}{1}{10}\claw{105}{1}{10}
\put(-10,-18){$\tinyM p_k+m_1$}\put(38,-18){$\tinyM p_{k-1}$}\put(70,-10){$\dots$}
\put(102,-18){$\tinyM  p_2$}
\put(-15,5){$\tinyM[m_kE_k]$}
\put(70,5){$\dots$}
\put(20,5){$\tinyM [m_{k-1}E_{k-1}]$}
\put(90,5){$\tinyM[m_2E_2]$}\put(125,-3){$\equiv$}

\put(142.5,-2.5){$\bullet$}\put(135,5){$\tinyM [m_2E_1]$}\put(165,-3){$+$}
\end{picture}
\begin{picture}(0,0)(-10,-3)
\put(-40,-3){$\suml^{k-1}_{i=2}\!\Big($}
\curve(0,0,40,0)\put(-2.5,-2.5){$\bullet$}\put(37.5,-2.5){$\bullet$}
\claw{0}{1}{10}\claw{40}{1}{10}
\put(-25,-20){$\tinyM m_i-m_1-\suml_{j>i}p_j$}\put(30,-17){$\tinyM m_1+\suml_{j>i}p_j$}
\put(-20,7){$\tinyM [m_iE_1]$}\put(25,7){$\tinyM [m_{i+1}E_2]$}\put(60,-3){$-$}
\end{picture}
\begin{picture}(0,0)(-100,-3)
\put(-2.5,-2.5){$\bullet$}\claw{0}{1}{10}\put(-15,7){$\tinyM [m_{i+1}E_1]$}
\put(20,-3){$\Big)$}
\end{picture}
\eeq

\

\

Recall that $m_k=km_1-p_1$, thus $m_k-jm_1=(k-j)m_1-p_1$.
This gives the statement.
\epr
\bex
1. If one applies part one of the theorem to the exceptional divisor of the first blowup of $(\C^2,0)$ then one gets equation (\ref{Eq.Spectral.Relation.Simplest.Case}).

2. The theorem gives also the useful recursive formula:\\
\beqm
\begin{picture}(0,0)(-20,-3)
\curve(0,0,125,0)
\put(-2.5,-2.5){$\bullet$}\put(37.5,-2.5){$\bullet$}\put(77.5,-2.5){$\bullet$}\put(122.5,-2.5){$\bullet$}
\claw{0}{1}{10}\claw{40}{1}{10}\claw{80}{1}{10}\claw{125}{1}{10}
\put(-2,-15){$\tinyM p_1$}\put(38,-15){$\tinyM p_k$}\put(75,-15){$\tinyM p_{k-1}$}\put(123,-15){$\tinyM p_2$}
\put(-10,5){$\tinyM[m_1E_1]$}\put(30,5){$\tinyM[m_kE_k]$}\put(65,5){$\tinyM[m_{k-1}E_{k-1}]$}
\put(115,5){$\tinyM[m_2E_2]$}\put(140,-3){$\equiv$}
\end{picture}
\begin{picture}(0,0)(-180,-3)
\curve(0,0,125,0)
\put(-2.5,-2.5){$\bullet$}\put(37.5,-2.5){$\bullet$}\put(77.5,-2.5){$\bullet$}\put(122.5,-2.5){$\bullet$}
\claw{0}{1}{10}\claw{40}{1}{10}\claw{80}{1}{10}\claw{125}{1}{10}
\put(-2,-15){$\tinyM p_1$}\put(38,-15){$\tinyM p_{k-1}$}\put(75,-15){$\tinyM p_{k-2}$}\put(123,-15){$\tinyM p_2$}
\put(-10,5){$\tinyM[m_1E_1]$}\put(20,5){$\tinyM[m_{k-1}E_{k-1}]$}\put(65,5){$\tinyM[m_{k-2}E_{k-2}]$}
\put(115,5){$\tinyM[m_2E_2]$}\put(150,-3){$+$}
\end{picture}
\\\\
\begin{picture}(0,0)(300,-3)
\curve(0,0,40,0)\put(-2.5,-2.5){$\bullet$}\put(37.5,-2.5){$\bullet$}
\claw{0}{1}{10}\claw{40}{1}{10}
\put(-20,7){$\tinyM [m_{k-1}E_1]$}\put(25,7){$\tinyM [m_kE_2]$}\put(60,-3){$-$}
\end{picture}
\begin{picture}(0,0)(210,-3)
\put(-2.5,-2.5){$\bullet$}\claw{0}{1}{10}\put(-15,7){$\tinyM [m_{k-1}E_1]$}\put(30,-3){$+$}
\end{picture}
\begin{picture}(0,0)(140,-3)
\put(-2.5,-2.5){$\bullet$}\claw{0}{1}{10}\put(-20,7){$\tinyM [(m_k-m_1)E_1]$}\put(30,-3){$-$}
\end{picture}
\begin{picture}(0,0)(70,-3)
\curve(0,0,40,0)\put(-2.5,-2.5){$\bullet$}\put(37.5,-2.5){$\bullet$}
\claw{0}{1}{10}\claw{40}{1}{10}
\put(-25,7){$\tinyM [(m_k-m_1)E_1]$}\put(30,7){$\tinyM [m_kE_2]$}
\end{picture}
\end{multline}
\eex

\subsection{Independence of the basic types}\label{Sec.Independence.Basic.Types}
Theorem \ref{Thm.Spectrum.Cp.Formulas} expresses the spectrum of $(C,0)$ through those of the "basic types" $C_{p,q}$. These basic types are linearly independent:

{\em There are non non-trivial identities of the form $\suml_{(p,q)}a_{p,q} Spec(C_{p,q})=0$ among the spectra of the basic types, $\{C_{p,q}\}_{\substack{p=0,1,\dots\\q=0,2,3,\dots}}$.}

We will prove this statement in the next versions of the paper.
Below we indicate the idea of the proof.

Consider some irreducible linear relation among the basic types. (Irreducible means that it does not decompose into the sum of relations involving fewer number of types.)

For some positive $m$ and $k$ such that $gcd(k,m)=1$. Consider the types that can contribute to $t^{\frac{k}{m}}$. These are:
\li the ordinary multiple point, $C_{m,0}$, it contributes $t^{\frac{k}{m}}(m-k-1)$.
Here the maximal possible spectral value is $\al_{max}=1-\frac{2}{m}$.
\li the types $C_{p,q}$ for $p+q=m$, they contribute $t^{\frac{k}{m}}(p-1-\lfloor\frac{pk}{m}\rfloor)$.
Here, for $p>0$, the maximal possible spectral value is  $\al_{max}\le 1-\frac{2}{m}$.
\li the types $C_{p,q}$ for $p+2q=m$, they contribute $t^{\frac{k}{m}}(q-\lceil\frac{qk}{m}\rceil)$.
 Here the maximal possible spectral value is $\al_{max}=1-\frac{3}{m}$.
\li the "external" types $C_{p,q}$, with either $(p+q)=ml$ or $(p+2q)=ml$, $l>1$.

Let $M$ be the maximal among the denominators that appear in this relation, i.e. the maximal among all the $(p+q)$'s and $(p+2q)$'s. Suppose $M\ge m>\frac{M}{2}$, then only the first three types above contribute to $t^{\frac{k}{m}}$ for $(k,m)=1$.

Therefore, for $\suml_{q>1} a_i C_{m-2q,q}+a_0C_{m,0}+\suml_{p<m-1} b_pC_{p,m-p}$, we get the vanishing condition on the coefficient of $t^{\frac{k}{m}}$:
\beq
\suml^{\lceil\frac{m}{2}\rceil-1}_{q=2}a_q(q-\lceil\frac{qk}{m}\rceil)+
\underbrace{a_{\frac{m}{2}}(\frac{m}{2}-1-\lfloor\frac{k}{2}\rfloor)}_{for\ m\in2\Z}+
a_0(m-k-1)+\suml^{m-2}_{p=2m-M} b_p(p-1-\lfloor\frac{pk}{m}\rfloor)=0,\quad k=1,\dots,m-1,\quad (m,k)=1
\eeq
By renaming $k\to m-k$ we get the condition on the coefficient of $t^{1-\frac{k}{m}}$:
\beq
\suml^{\lceil\frac{m}{2}\rceil-1}_{q=2}a_q\lfloor\frac{qk}{m}\rfloor+
\underbrace{a_{\frac{m}{2}}(\lceil\frac{k}{2}\rceil-1)}_{for\ m\in2\Z}+
a_0(k-1)+\suml^{m-2}_{p=2m-M} b_p(\lceil\frac{pk}{m}\rceil-1)=0,\quad k=1,\dots,m-1,\quad (m,k)=1
\eeq
Suppose $k$ divides $m$ but is prime. Then $t^{1-\frac{k}{m}}=t^{1-\frac{1}{m'}}$, for some $m'\in\Z$. But such a power is not realized as $\al_{max}$ for any other types. Thus, in this case the sum of coefficients of
$t^{1-\frac{k}{m}}$ must vanish as well, i.e. the last equation must hold for $k$ satisfying: either
$(m,k)=1$ or $m\in k\Z$.

For $k=1$ the last equation is trivial. For $m=M$ or $m=M-1$ the terms $b_p$ do not participate. Thus, for $m=M$, $(M-1)$, by taking $k=2$ we get: $a_0=0$.

Thus, for $m=M$ or $m=M-1$ the variables are $\{a_q\}_{q=2,\dots,\lceil\frac{m}{2}\rceil-1}$. The number of equations can be computed using Euler's function, $\phi(m)=m\prodl_{m_i|m}(1-\frac{1}{m_i})$. To this
 one adds the number of prime divisors of $m$.

In many cases (e.g. $m$ is prime or $m$ has no small enough prime factors) the total $\phi(m)+r$ is bigger than the number of variables. The system of equations which we get is trapezoidal.

\subsection{The rings \Types and \Spec}
Consider the collection of all the (embedded topological) singularity types $S_i$ of plane curve germs $(C,0)\sset(\C^2,0)$.
Denote by \Types the abelian group freely generated by the formal linear combinations, $\sum a_i S_i$.
 This abelian group admits the natural "merging" product, $(S_1,S_2)\to S_1\otimes S_2$. Take the generic enough representatives, $(C_i,0)$ of $S_i$, in particular they have no common tangents.
 Define $S_1\otimes S_2$ to be the singularity type of $(C_1\cup C_2,0)$.
 The singularity type of this germ is fully determined by $S_1,S_2$, i.e. does not depend on the choice of the representatives $(C_i,0)$.

Consider the collection of all the spectra of (reduced) plane curve singularities. Denote by \Spec the abelian group of all the (formal finite) linear combinations $\sum a_i Spec(C_i,0)$. Using theorem \ref{Thm.Spectrum.Cp.Formulas}  and assuming independence of basic types (section \ref{Sec.Independence.Basic.Types}) we get:
\bcor
 \Spec is the quotients of (the abelian group) \Types by all the recombination relations. \Spec is freely generated by the basic types $\{C_{p,q}\}$.
\ecor
When does the product of \Types descend to a product on \Spec? To descend it must preserve the subgroup of the relations.
\bex
1. Consider the relation 1 of theorem \ref{Thm.Spectrum.Cp.Formulas}, multiply it by a smooth branch. We get:\\
\beq
\begin{picture}(0,0)(180,0)
\curve(0,0,40,0)\put(20,0){$\curlyvee$}\put(30,0){$\curlyvee$}
\put(20,8){$\overbrace{}$}  \put(25,17){$\tinyM D$}
\put(-5,-3){\blackboxh}\put(-28,-2){$\tinyM m'_i E_i$}
\put(50,0){$-$}
\end{picture}
\begin{picture}(0,0)(80,0)
\curve(0,0,40,0)\freev{30}{-5}{10}
\put(20,8){$\overbrace{}$}  \put(15,17){$\tinyM deg(D\cap E_i)$}
\put(-5,-3){\blackboxh}\put(-28,-2){$\tinyM m'_i E_i$}
\put(50,0){$\equiv??\equiv$}
\end{picture}
\begin{picture}(0,0)(-40,0)
\curve(0,0,40,0)\put(20,0){$\curlyvee$}\put(30,0){$\curlyvee$}
\put(20,8){$\overbrace{}$}  \put(25,17){$\tinyM D$}
\freev{10}{-5}{10}\put(-30,-10){$\tinyM (m_i+1) E_i$}
\put(50,0){$-$}
\end{picture}
\begin{picture}(0,0)(-140,0)
\curve(0,0,40,0)\freev{10}{-5}{10}\freev{30}{-5}{10}
\put(5,8){$\overbrace{\quad\quad\quad}$}  \put(15,17){$\tinyM (m_i+1)$}
\put(-30,-10){$\tinyM (m_i+1) E_i$}
\end{picture}
\eeq\\
Suppose in the "preliminary history", i.e. $\blackboxh$, we had to blowup at least once at the intersection point of two exceptional divisors. (This happens when there is a non-smooth branch.) Then $m'_i>m_i+1$ and the formula does not hold anymore.

2. Let $\Spec^{s.b.}\sset\Spec$ denote the abelian subgroup generated by $Spec(C,0)$ for curve singularities with smooth branches.  Then $\Spec^{s.b.}$ is the quotient of $\Types^{s.b.}$ by the relations as above. For smooth branches one has $m'_i=m_i+1$, i.e. the relations admit multiplication by the smooth branches.
\eex
In view of these examples, one should either restrict to $\Spec^{s.b.}$ or to define the product in a different way.
 We define the product directly on $\Spec$.

 First define it on the generators: $C_{p',0}\otimes C_{p,q}=C_{p+p',q}$ and
 $C_{p,q}\otimes C_{p',q'}=C_{p+p'+q',q}+C_{p+p'+q,q'}-C_{p+p'+q'+q,0}$.
  (The later relation is forced by the compatibility with recombination formulae.)
Now extend the product by linearity to \Spec. Add to \Spec (formally) the unit element, $\one$.
\bcor
With this product $(\Spec,\one)$ becomes an associative commutative ring over $\Z$, generated by the types
$C_{1,0}$ and $\{C_{0,q}\}_{q\ge2}$. For the types with smooth branches the product is "geometric", i.e. descends from $\Types^{s.b.}$.
\ecor
Note that \Spec is not a polynomial ring, its the generators are algebraically related:
\beq
C_{0,q}\otimes C_{0,q'}=C_{1,0}^{\otimes q')}\otimes C_{0,q}+C_{1,0}^{\otimes q)}\otimes C_{0,q'}-C_{1,0}^{\otimes(q+q')}.
\eeq
\bpr
We should prove commutativity, distributivity, associativity. As the product is defined by generators,
 it is enough to check commutativity and associativity for generators only.
 The only non-trivial cases are the triple product, $C_{p,0}\otimes C_{p',q'}\otimes C_{p'',q''}$ and
  $C_{0,q}\otimes C_{0,q'}\otimes C_{0,q''}$, they are checked directly.

Finally, the distributivity holds as the basic types are linearly independent.
\epr

\subsection{The non-additivity of spectral pairs}\label{Sec.Spectral.Pairs.Non.additive}
The spectral pairs  are more delicate invariants, they are not additive, even in the simplest cases.
\\\\
\parbox{12cm}
{For example, consider the singularity of the type $(C,0)=\{\prodl^r_{i=1}(l_i^p+y^{2p}+x^{2p})=0\}\sset(\C^2,0)$, where $\{l_i\}$ are pairwise non-proportional linear forms. The resolution is on the diagram, the spectral pairs are:}
\begin{picture}(0,0)(-30,0)
\linethickness{0.8mm}\curve(0,0,120,0)\curve(35,-5,35,30)\curve(75,-5,75,30)\curve(105,-5,105,30)\thinlines
\freev{15}{-5}{10}\freeh{30}{15}{10}\freeh{70}{15}{10}\freeh{100}{15}{10}
\put(-25,0){$\tinyM prE_1$}\put(20,-15){$\tinyM(pr+r)E_2$}\put(90,-15){$\tinyM(pr+r)E_{r+1}$}
\put(50,15){$\dots$}
\end{picture}
\beq
Spp(C,0)=r\suml_{0<s<pr+p}\Big(p-1-\lfloor\frac{s}{r+1}\rfloor\Big)\Big(t^{[-\frac{s}{pr+r},1]}+t^{[\frac{s}{pr+r},1]}\Big)
+(r-1)\suml_{0<s<p}\Big(t^{[-\frac{s}{p},2]}+t^{[\frac{s}{p},0]}\Big)+(pr-1)t^{[0,1]}
\eeq
The natural wish is to use the "tangential decomposition", as in equation (\ref{Eq.Spectral.Relation.Simplest.Case}):
\beq
Spp(C,0)=?= r Spp\Big((x^{p(r-1)}+y^{p(r-1)})(x^p+y^{2p})\Big)-(r-1)Spp(x^{pr}+y^{pr})
\eeq
But this decomposition does not hold for the spectral pairs. Indeed:
\beqm
Spp\Big((x^{p(r-1)}+y^{p(r-1)})(x^p+y^{2p})\Big)=
\suml_{0<s<pr+p}\Big(p-1-\lfloor\frac{s}{r+1}\rfloor\Big)\Big(t^{[-\frac{s}{pr+r},1]}+t^{[\frac{s}{pr+r},1]}\Big)
+\\+\suml_{0<s<pr}\Big(pr-p-1-s+\lceil\frac{s}{q}\rceil\Big)\Big(t^{[-\frac{s}{pr},1]}+t^{[\frac{s}{pr},1]}\Big)
+\suml_{0<s<p}\Big(t^{[-\frac{s}{p},2]}+t^{[\frac{s}{p},0]}\Big)+(pr-1)t^{[0,1]}
\end{multline}
while
\beq
Spp(x^{pr}+y^{pr})=\suml_{0<s<pr}(pr-1-s)\Big(t^{[-\frac{s}{pr},1]}+t^{[\frac{s}{pr},1]}\Big)+(pr-1)t^{[0,1]}
\eeq
Thus
$Spp(C,0)\neq r Spp\Big((x^{p(r-1)}+y^{p(r-1)})(x^p+y^{2p})\Big)-(r-1)Spp(x^{pr}+y^{pr})$

\subsection{The non-additivity of the Seifert form}
The Seifert form is not additive either. Even a weaker invariant: the collection of equivariant signatures of the Seifert form, \cite[Proposition 3.1]{Schr.Stee.Stev.} is not additive.

In the example of section
\ref{Sec.Spectral.Pairs.Non.additive} one has:
\beq\ber
for\ \prodl^r_{i=1}(l_i^p+y^{2p}+x^{2p}):\quad \si^-_{\frac{s}{pr+p}}=r\Big(p-\lfloor\frac{s}{r+1}\rfloor-\lfloor-\frac{s}{r+1}\rfloor\Big)
\\
for\ (x^{p(r-1)}+y^{p(r-1)})(x^p+y^{2p}):\quad \si^-_{\frac{s}{pr+p}}=p-\lfloor\frac{s}{r+1}\rfloor-\lfloor-\frac{s}{r+1}\rfloor,\quad
\si^-_{\frac{s}{pr}}=pr-p-2s+\lfloor\frac{s}{q}\rfloor-\lfloor-\frac{s}{q}\rfloor
\\
for\ (x^{pr}+y^{pr}):\quad \si^-_{\frac{s}{pr}}=pr-2s
\eer\eeq
Thus $\si^-_{\frac{s}{pr+p}}\Big(\prodl^r_{i=1}(l_i^p+y^{2p}+x^{2p})\Big)=
r\si^-_{\frac{s}{pr+p}}\Big((x^{p(r-1)}+y^{p(r-1)})(x^p+y^{2p})\Big)$,
\\ but
 $0\neq r\si^-_{\frac{s}{pr}}\Big((x^{p(r-1)}+y^{p(r-1)})(x^p+y^{2p})\Big)-(r-1)\si^-_{\frac{s}{pr}}\Big(x^{pr}+y^{pr}\Big)$.

\section{Applications}

\subsection{Singularity types with coinciding spectrum}\label{Sec.Sings.With.Coinciding.Spectra}
Theorem \ref{Thm.Spectrum.Cp.Formulas} gives a very simple way to construct curve singularities of different topological types, whose spectra coincide. For example:
\beq
\begin{picture}(0,0)(130,40)
\linethickness{0.8mm}\curve(0,0,60,0)\curve(10,-5,10,45)\curve(50,-5,50,45)
\thinlines\freeh{5}{15}{10}\freeh{5}{40}{10}
\put(50,10){$\prec$}\put(50,35){$\prec$}
\put(-25,0){$\tinyM mE_1$}\put(20,-15){$\tinyM E_2$}
\put(90,10){$\equiv$}
\end{picture}
\begin{picture}(0,0)(0,40)
\linethickness{0.8mm}\curve(0,0,60,0)\curve(10,-5,10,45)\curve(50,-5,50,45)
\thinlines\freeh{5}{10}{10}\freeh{50}{40}{10}
\put(10,35){$\prec$}\put(50,10){$\prec$}
\put(-25,0){$\tinyM mE_1$}\put(20,-15){$\tinyM E_2$}
\end{picture}
\eeq

\vspace{2cm}
though the topological types are different.

\subsection{What is determined by spectrum?}\label{Sec.What.Is.Determined.By.Spectrum}
Recall that even the spectral pairs do not determine the (embedded topological) singularity types, \cite{Schr.Stee.Stev.}. On the other hand,
 the Seifert form (which is determined by the spectral pairs) determines the intersection multiplicities of the branches, \cite{Kaenders1996}.

Using theorem \ref{Thm.Spectrum.Cp.Formulas}  and linear independence of the basic types we get:
\bcor
1. The multiplicity of $(C,0)$ is determined by $Spec(C,0)$.
\\2. If $(C,0)$ has only smooth branches or $(C,0)$ is of the 'intermediate type', $C_{p_1,\dots,p_k}$, then the multiplicities of all the exceptional divisors of the minimal (good) resolution of $(C,0)$ are determined by $Spec(C,0)$.
\ecor
\bpr
1. Expand $Spec(C,0)$ into the linear combination of the basic types, as in theorem \ref{Thm.Spectrum.Cp.Formulas}. To restore the multiplicity  it is enough to prove that it does not cancel out.

First assume that the tangent cone  $T_{(C,0)}$ has at most one multiple line (and possibly some other reduced lines). Then in the expansion one
 does not use the recombinations on $E_1$. Thus, in the resulting collection of the basic types the edge of multiplicity $mult(C,0)$ occurs precisely one, with coefficient $+1$.
 So it cannot cancel with anything else.

In the general case, i.e. several multiple lines in $T_{(C,0)}$, first use the recombination on $E_1$ to get the formula as in the intro, equation (\ref{Eq.Spectral.Relation.Simplest.Case}).
Now expand each $Spec(C^{(i)},0)$ further. Then the resulting expansion contains $\{C_{p,q}\}_{\substack{q=2,\dots,\\p+q=m}}$
 with positive coefficients and possibly $C_{m,0}$ with negative coefficient. Thus the terms that involve $m$ cannot cancel.

2. If $(C,0)$ has only smooth branches then, as one sees from part 3 of theorem \ref{Thm.Spectrum.Cp.Formulas}, no cancelation is possible. Namely, in the decomposition $(C,0)=\suml_{p_1,p_2} (C_{p_1,p_2},0)-\suml_p (C_{p,0},0)$ all the participating summands are determined by $Spec(C,0)$.

If $(C,0)$ is of the type $C_{p_1,\dots,p_k}$ then use the expansion of part 4 of theorem \ref{Thm.Spectrum.Cp.Formulas}. The right hand side of equation (\ref{Eq.Spec.Decomposition.of.Intermediate.Types}) contains
\beq
\begin{picture}(0,0)(200,-3)
\put(-23,-3){$\suml_i$}
\curve(0,0,40,0)\put(-2.5,-2.5){$\bullet$}\put(37.5,-2.5){$\bullet$}
\claw{0}{1}{10}\claw{40}{1}{10}
\put(-24,-20){$\tinyM m_i-m_1-\suml_{j>i}p_i$}\put(30,-17){$\tinyM m_1+\suml_{j>i}p_i$}
\put(-15,7){$\tinyM [m_iE_1]$}\put(25,7){$\tinyM [m_{i+1}E_2]$}
\end{picture}
\begin{picture}(0,0)(50,-3)
\curve(0,0,40,0)
\put(-2.5,-2.5){$\bullet$}\put(37.5,-2.5){$\bullet$}
\claw{0}{1}{10}\claw{40}{1}{10}
\put(-2,-15){$\tinyM p_1$}\put(25,-15){$\tinyM m_1-p_1$}
\put(-10,5){$\tinyM [m_1E_1]$}\put(20,5){$\tinyM [(2m_1-p_1)E_2]$}
\end{picture}
\begin{picture}(0,0)(-100,-3)
\put(-30,-3){$\suml_i$}
\curve(0,0,70,0)\put(-2.5,-2.5){$\bullet$}\put(67.5,-2.5){$\bullet$}
\claw{0}{1}{10}\claw{70}{1}{10}
\put(-14,-17){$\tinyM (i-1)m_1-p_1$}\put(65,-17){$\tinyM m_1$}
\put(-20,7){$\tinyM [(im_1-p_1)E_1]$}\put(35,7){$\tinyM [(i+1)m_1-p_1)E_2]$}
\end{picture}
\eeq

\

\

We claim that these summands cannot cancel. Again, as these are basic types (linearly independent),
the cancelation can occur only if some multiplicities coincide, i.e.
\beq
\bpm m_i=m_1\\m_{i+1}=2m_1-p_1\epm\quad \rm{or}\ \
\bpm m_i=jm_1-p_1\\m_{i+1}=(j+1)m_1-p_1\epm\quad \rm{or}\ \
\bpm jm_1-p_1=m_1\\(j+1)m_1-p_1=2m_1-p_1\epm
\eeq
Recall that for $i\ge1$ $m_{i+1}=m_i+\suml_{j>i}p_j$ and $m_1>\suml_{j\ge1}p_j$.
Thus we get contradiction in each of the cases.
\epr

\subsection{A relative point of view}\label{Sec.Relative.Point.of.View}
The relations  of theorem \ref{Thm.Spectrum.Cp.Formulas} suggest that in many cases it is useful to consider not the  singularity itself, but the difference: $[(C,0)]:=(C,0)-(x^m=y^m)$, where $m=mult(C,0)$.
 Then, for any invariant that is compatible with the \cp\ formulae we define $inv[(C,0)]:=inv(C,0)-inv(x^m=y^m)$.
Thus, e.g. to prove an inequality $inv(C,0)>0$ it is enough to verify: $inv(x^m-y^m)>0$ and $inv[(C,0)]>0$.
 In many cases it is simpler to check  $inv[(C,0)]>0$ than $inv(C,0)>0$.

\subsection{Stabilizations of curve singularities with small $|\mu_0+\mu_+|$}\label{Sec.Surfaces.Sings.low.muPlusZero}
Consider an isolated surface singularity $(X,0)\sset(\C^3,0)$. The intersection form on the homology of the Milnor fibre, $H_*(X_F,\R)$, is determined by the triple
 $(\mu_-,\mu_0,\mu_+)$, the number of negative/ zero/positive indices.
An old question is the classification of surface singularities with small values of $\mu_0+\mu_+$. The cases $\mu_0+\mu_+\le 2$ have been classified by Ebeling and others.

Recall that for $(X,0)\sset(\C^3,0)$: $\mu_-=|Spec(X,0)\cap(0,1)|$ while $\mu_++\mu_0=2|Spec(X,0)\cap(-1,0]|$.
Let $m$ be the multiplicity of $(X,0)$, by the semicontinuity of the spectrum we have $\mu_-(X,0)\ge\mu_-(x^m+y^m+z^m=0)$ and similarly for $(\mu_++\mu_0)$. Thus, for $[(X,0)]=(X,0)-(x^m+y^m+z^m=0)$ one has $\mu_-[(X,0)]\ge0$ and $(\mu_++\mu_0)[(X,0)]\ge0$, and one should study the cases of low values of
 $(\mu_++\mu_0)[(X,0)]$.

We consider a particular case: surface singularities which are stabilizations of curve singularities, $(X,0)=\{f(x,y)+z^2=0\}$.
 By the Thom-Sebastiani formula, if $Spec(C,0)=\sum t^{\al_i}$ then $Spec(X,0)=\sum t^{\al_i+\frac{1}{2}}$ .
 Thus
  $(\mu_++\mu_0)[(X,0)]=2|Spec[(C,0)]\cap(-1,-\frac{1}{2}]|$.
 Now, the first statement of theorem \ref{Thm.Spectrum.Cp.Formulas} reduces the study of  $|Spec[(C,0)]\cap(-1,-\frac{1}{2}]|$ to the "intermediate types", $(C_{p_1,\dots,p_k},0)$.
In particular we record:
\bel
Let $(C_{p_1,p_2},0)$ be the plane curve singularity of the type $\{(x^{p_1}+y^{p_1})(x^{p_2}+y^{2p_2})=0\}$.
 Then
 \[
 (\mu_++\mu_0)[(C_{p_1,p_2},0)]=\bin{\lceil\frac{p_2}{2}\rceil}{2}-\frac{1-(-1)^{p_1p_2}}{2}\frac{p_2-1}{2}.\]
\eel

\subsection{Strengthening of Givental's bound}\label{Sec.Givental.Bound}
For a plane curve singularity consider the positive part of the spectrum:
$Spec(C,0)\cap(0,1)=\{0<\al\le\al_2\le\cdots\le\al_{k}<1\}$.
In \cite{Givental} one proves: $\al_{i+1}+\al_{k-i}<1$ for all $0<i<k$.
We strengthen this bound as follows.
\bprop
Let $Spec(C,0)=\suml_{i>0} a_i t^{\pm \al_i}+a_0t^0$, where $0<\al_1\lneqq\al_2\lneqq\cdots\lneqq\al_k$.
 Consider the (minimal embedded good) resolution of $(C,0)$, with the exceptional divisors $\{m_iE_i\}$.
 Denote by $r$ the maximal number of distinct multiplicities $\{m_i\}_{i\in I}$
 that satisfy: $\{\frac{1}{m_i}<\frac{2}{m_j}\}_{\substack{\forall i,j\in I\\i\neq j}}$.
  Then $\al_{i+r+1}+\al_{k-i}\le1$, and the equality holds iff $(C,0)$ is an ordinary multiple point.
\eprop
Note that the number $r$ of the proposition is always at least 1, in most cases it is bigger than 1.
\\When compared to the initial Givental's bound, the significant strengthening is in the condition $\al_i\lneqq\al_{i+1}$. (Recall that usually the spectral values that are close to 0 come with large multiplicities, while those close to 1 usually have multiplicity one.)
\\\bpr
First we check the bound for the \omp, cf. example \ref{Ex.OMP.Basic.Types.Spectrum}. Here $\al_i=\frac{i}{m}$ for $i=1,\dots,m-2$, thus  $\al_{i+2}+\al_i=1$ for $i\ge0$.

{\bf Step 1.} We start with elementary remarks before addressing the general curve singularity.
 Note that the bound $\al_{i+r+1}+\al_{k-i}\le1$ does not involve $\al_1,\dots,\al_r$. Thus we often ignore/erase  these first few terms. We say that the sequence $\al_1<\al_2<\cdots<\al_k$ satisfies $\cG$-condition if $\al_{i+1}+\al_{k-i}<1$ for $0\le i<\frac{k-1}{2}$ and if $k$ is odd then $\al_{\frac{k-1}{2}}\le\frac{1}{2}$.

Suppose some sequence $0<\al_1<\cdots<\al_k$ satisfies $\cG$ and $0<\be<\al_1$. Then
$0<\be<\al_1<\cdots<\al_k$ satisfies $\cG$ as well. Further, $0<\al_1<\cdots<\al_k$ satisfies $\cG$ iff
 $\al_1+\al_k<1$ and $\al_2<\cdots<\al_{k-1}$ satisfies $\cG$.

Given two  increasing sequences, $\{\al_i\}_{i=1,\dots,k}$ and $\{\be_i\}_{i=1,\dots,s}$, consider their union (with repetitions omitted):
 $0<\ga_1<\cdots<\ga_t$. We claim that if $\{\al_i\}$, $\{\be_i\}$ satisfy $\cG$ then $\{\ga_i\}$ satisfies $\cG$ as well.
  The proof goes by reduction in the length of the sequence, note that the statement is trivial if one of the sequences is one of the sequences is empty.
When merging the two sequences into $\{\ga_i\}$ two cases are possible:
\li $\ga_1$, $\ga_t$ come from the same sequence, e.g. $\ga_1=\al_1$ and $\ga_t=\al_k$. Then $\ga_1+\ga_t<1$
 and $\ga_1<\cdots<\ga_t$ satisfies $\cG$ iff $\ga_2<\cdots<\ga_{t-1}$ satisfies $\cG$. Thus one erases
 $\al_1,\al_k$ from $\{\al_i\}$ and repeats the procedure, merging $\al_2,\dots,\al_{k-1}$ and $\{\be_j\}$.
\li
$\ga_1$, $\ga_t$ do not come from the same sequence, e.g. $\ga_1=\al_1$ and $\ga_t=\be_s$. Thus $\al_1<\be_1$ and $\al_k<\be_s$, so $\be_1+\al_k<\be_1+\be_s<1$ and $\al_1+\be_s<\be_1+\be_s<1$. Then, instead of $\{\al_i\}$ and $\{\be_i\}$ one can consider $\al_1<\al_2<\cdots<\al_{k-1}<\be_s$ and $\be_1,\be_2,\cdots,\be_{s-1},\al_k$ (the last sequence is reordered if needed).
As $\al_1+\be_s<1$ the sequence $\al_1<\al_2<\cdots<\al_{k-1}<\be_s$ satisfies $\cG$. Similarly, as $\be_1+\al_k<1$, the second sequence satisfies $\cG$. Thus $\{\ga_i\}$ is the result of merging of these two new sequences, and now
 $\ga_1$, $\ga_t$ come from the same sequence ($\al_1<\al_2<\cdots<\al_{k-1}<\be_s$). So, we are in the first case and proceed as above (erase $\ga_1$, $\ga_t$ and repeat the process).

{\bf Step 2.} Consider the spectrum $Spec(C,0)=\suml_i\suml_k d_{k,i}t^{\frac{k}{m_i}}$, cf. lemma \ref{Thm.Spectrum.Additive.Formula}. As one sees from the formula: $d_{k,i}\ge0$ except possibly for one case, when the first exceptional divisor $E_i$ intersects no branches of the curve (i.e. $p_1=0$) and only one among the other exceptional divisors (i.e. $E_{i>1}\cap E_1=\empty$ except for one $E_i$).

Suppose $d_{k,i}\ge0$ then for each $E_i$ we get $\suml_k d_{k,i}t^{\frac{k}{m_i}}$,
 which produces the sequence $0<\frac{1}{m_i}<\cdots<\frac{k_{max}}{m_i}$. As $(C,0)$ is not an \omp, $k_{max}<m_i-2$, thus the cropped sequence $\frac{2}{m_i}<\cdots<\frac{k_{max}}{m_i}$ satisfies $\cG$. Go over all $E_i$ and merge  the corresponding cropped sequences, the resulting $\ga_1<\cdots<\ga_k$ satisfies $\cG$ again.

 Among the values $\{\frac{1}{m_i}\}_i$ choose those that are not smaller than $\ga_1$ and add them to the sequence $\{\ga_i\}$. As all of  them are not bigger than $\frac{1}{2}$ the obtained sequence satisfies $\cG$.
Finally, by the assumption of the proposition, at least $r$ values among $\{\frac{1}{m_i}\}_i$ are distinct and smaller than $\ga_1$. Thus, after they are added, we get the sequence $0<\{\frac{1}{m_i}\}_i<\ga_1<\cdots<\ga_k$ that satisfies the claimed statement.

When $d_{k,1}=-1$, for $E_1$, this just decreases the multiplicities of spectral values on the corresponding $E_i$. As the spectral values remain the same, this does not break the $\cG$ condition.
\epr

\section{The generalized Durfee bound}\label{Sec.Durfee.Bound.Generalized}
In \cite{Durfee1978} the bound was stated on the singularity genus and the Milnor number of  normal surface singularities in $(\C^3,0)$: $\mu\ge 6p_g$.
By now the general conjecture for isolated hypersurface singularity $(X,0)\sset(\C^n,0)$ is: $\mu(X,0)\ge n!p_g(X,0)$.  It has been verified in several cases, see \cite{Kerner.Nemethi.2014}.

Recall that the singularity genus is expressible from the spectrum, $p_g(X,0)=|Spec(X,0)\cap(-1,0]|$, \cite{Saito1981}. Thus Durfee conjecture can be stated as: $\mu(X,0)\ge n!|Spec(X,0)\cap(-1,0]|$.
We propose the more general
\bconj
Let $-1<\al\le0$ and $n\ge2$, then $\frac{\mu(X,0)}{n!}\ge\frac{|Spec(X,0)\cap(-1,-\al]|-C_n}{(1-\al)^n}$, the constant $C_n$ depends on $n$ only.
\econj
We verify this conjecture in several cases.

\subsection{\Nnd\ hypersurface singularities with large enough Newton diagrams}
In \cite{Kerner.Nemethi.2014} we prove the bound $\mu\ge n!p_g$ for large enough diagrams. More precisely:
 for any \ND\ $\Ga\sset\R^n_{\ge0}$ and any $t\gg0$ the bound holds for \Nnd\  singularities with the diagram $t\Ga$.
Using this result we can establish the generalized Durfee bound for this class of singularities.
\bprop
Fix some $0\le \al<1$, suppose the diagram $(1-\al)\Ga$ is large enough. Then
 $\frac{\mu(X,0)}{n!}>\frac{\Big|Spec(X,0)\cap(-1,-\al]\Big|}{(1-\al)^n}$.
\eprop
\bpr
We use the expression for the spectrum in terms of $\Ga\sset\R^n_{\ge0}$. Then the values of $Spec(X,0)\cap(-1,-\al]$ correspond precisely to the points of $\Z^n_{>0}$, which are on or below the scaled diagram $(1-\al)\Ga$. Therefore
$|Spec(X,0)\cap(-1,-\al]|=p_g((1-\al)\Ga)$ and by \cite{Kerner.Nemethi.2014}, for $\Ga$ large enough,
$n!p_g((1-\al)\Ga)\le\mu((1-\al)\Ga)$. Thus, using Kouchnirenko's formula:
\beq
n!|Spec(X,0)\cap(-1,-\al]|\le \mu((1-\al)\Ga)=\suml^n_{i=0}(-1)^i(n-i)!Vol_{n-i}((1-\al)\Ga)=
(1-\al)^n\suml^n_{i=0}\frac{(-1)^i(n-i)!Vol_{n-i}(\Ga)}{(1-\al)^i}
\eeq
Here $Vol_i(\Ga)$ is the total volume of the intersection of $\Ga$ with all the $i$'th dimensional coordinate subspaces of $\R^n$, in particular $Vol_0(\Ga)=1$.
Present the later term in the form $(1-\al)^n\Big(\mu(\Ga)+\suml^{n}_{i=1}(-1)^i(n-i)!Vol_{n-i}(\Ga)(-1+\frac{1}{(1-\al)^i})\Big)$ and note that
 for large enough $\Ga$: $Vol_{n-1}(\Ga)\gg Vol_{n-2}(\Ga)\gg\cdots\gg Vol_{0}(\Ga)$. Thus
 $\suml^{n}_{i=1}(-1)^i(n-i)!Vol_{n-i}(\Ga)\Big(-1+\frac{1}{(1-\al)^i}\Big)\le0$, hence
 $n!|Spec(X,0)\cap(-1,-\al]|\le (1-\al)^n\mu(\Ga)$.
\epr

\subsection{Plane curve singularities}
\bprop
$\mu(C,0)+\frac{2mult(C,0)}{(1-\al)^2}>\frac{2}{(1-\al)^2}|Spec(C,0)\cap(-1,-\al]| $
\eprop
This bound is far from being sharp, we hope to improve it in the next versions of the paper.
\bpr
As explained in section \ref{Sec.Relative.Point.of.View}, it is enough to check the bound for the ordinary multiple point and for $[C_{p_1,\dots,p_r}]=C_{p_1,\dots,p_r}-C_{p_1+\suml^r_{i>1}(i-1)p_i}$, where
$C_{p_1+\suml^r_{i>1}(i-1)p_i}$ denote the \omp,  $(y^{p_1+\suml^r_{i>1}(i-1)p_i}-x^{p_1+\suml^r_{i>1}(i-1)p_i})$,
while $C_{p_1,\dots,p_r}=\{(y^{p_1}+x^{p_1})\prodl^{r}_{i=2}(y^{(i-1)p_i}-x^{ip_i})=0\}\sset(C^2,0)$.

For the ordinary multiple point, $x^m-y^m$: $\mu=(m-1)^2$, $Spec=\suml_{|k|<m}t^{\frac{k}{m}}(m-|k|-1)$. Thus
\beq
|Spec\cap(-1,-\al]|=\suml^{m-1}_{k=\lfloor\al m\rfloor}(m-k-1)=\bin{m-\lfloor\al m\rfloor}{2}=\bin{m-\al m}{2}+
\frac{\{\al m\}(2m-2\al m-1)+\{\al m\}^2}{2}.
\eeq
Altogether:
\beq
\mu-\frac{2}{(1-\al)^2}|Spec\cap(-1,-\al]|=1+\frac{\{\al m\}(1-\{\al m\})}{(1-\al)^2}+m\frac{2\al-1-2\{\al m\}}{(1-\al)}>-\frac{2m}{(1-\al)}
\eeq

The formula for the spectrum is written in lemma \ref{Thm.Spectrum.Additive.Formula}. In our case it gives:
\beqm
Spec\Big(y^{p_1}-x^{p_1})\prodl^{r}_{i=2}(y^{(i-1)p_i}-x^{ip_i})\Big)=
\suml_{|k|<m_1}t^{\frac{k}{m_1}}\Big(p_1-1-\frac{p_1k}{m_1}+\{-\frac{m_rk}{m_1}\}\Big)
+\suml_{|k|<m_2}t^{\frac{k}{m_2}}\Big(p_2-\frac{p_2k}{m_2}-\{\frac{m_3k}{m_2}\}\Big)+
\\+\suml^r_{i=3}\suml_{|k|<m_i}t^{\frac{k}{m_i}}\Big(p_i-\frac{p_ik}{m_i}-\{\frac{m_{i+1}k}{m_2}\}
+\{-\frac{m_{i-1}k}{m_2}\}\Big),\quad here\ m_{r+1}:=m_1
\end{multline}
Therefore
\beqm
\Big|Spec\Big(y^{p_1}-x^{p_1})\prodl^{r}_{i=2}(y^{(i-1)p_i}-x^{ip_i})\Big)\cap(-1,-\alpha]\Big|=
\suml_{i\ge1}\suml^{m_i-1}_{k\ge\lfloor\al m_i\rfloor}\Big(p_i-\frac{p_ik}{m_i}-\underbrace{1}_{only\ for\ i=1}\Big)+\\
+\underbrace{\suml^{m_1-1}_{k\ge\lfloor\al m_1\rfloor}\{-\frac{m_rk}{m_1}\}-
\suml^{m_2-1}_{k\ge\lfloor\al m_2\rfloor}\{\frac{m_3k}{m_2}\}
+\suml^r_{i=3}\suml^{m_i-1}_{k\ge\lfloor\al m_i\rfloor}\Big(\{-\frac{m_{i-1}k}{m_i}\}-\{\frac{m_{i+1}k}{m_i}\}\Big)}_{Frac^{(1)}_{m_1,\dots,m_r}}
\end{multline}
Now we can write the difference:
\beqm
\Big|Spec\Big(y^{p_1}-x^{p_1})\prodl^{r}_{i=2}(y^{(i-1)p_i}-x^{ip_i})\Big)\cap(-1,-\alpha]\Big|
-\Big|Spec\Big((y^{p_1+\suml^r_{i>1}(i-1)p_i}-x^{p_1+\suml^r_{i>1}(i-1)p_i})\Big)\cap(-1,-\alpha]\Big|=\\=
\suml_{i\ge2}\suml^{m_i-1}_{k\ge\lfloor\al m_i\rfloor}\Big(p_i-\frac{p_ik}{m_i}\Big)-
\suml^{m_1-1}_{k\ge\lfloor\al m_1\rfloor}\Big(m_1-p_1)(1-\frac{p_1k}{m_1}\Big)+Frac_{m_1,\dots,m_r}=\\=
\suml_{i\ge2}\frac{p_i}{m_i}\bin{m_i-\lfloor\al m_i\rfloor+1}{2}-
\frac{m_1-p_1}{m_1}\bin{m_1-\lfloor\al m_1\rfloor+1}{2}+Frac^{(1)}_{m_1,\dots,m_r}
\end{multline}
Present $\lfloor\al m_i\rfloor=\al m_i+\{\al m_i\}$ to get:
\beqm
\De(Spec)=(1-\al)^2\frac{\suml_{i\ge1}p_im_i-m^2_1}{2}+(1-\al)\frac{\suml_{i\ge1}p_i-m_1}{2}+
\\+\underbrace{\suml_{i\ge1}\Big(p_i\{\al m_i\}+p_i\frac{\{\al m_i\}+\{\al m_i\}^2}{2m_i}\Big)-
\Big(m_1\{\al m_1\}+\frac{\{\al m_1\}+\{\al m_1\}^2}{2}\Big)}_{Frac^{(2)}}+Frac^{(1)}_{m_1,\dots,m_r}
\end{multline}
Now use $m_{i+1}-m_i-m_1=\suml_{j\ge i+1}p_j$, $i>1$, and $m_2-m_1=\suml_{j\ge2}p_j$
to get: $m_i-im_1=\suml_{j\ge2}p_j+\suml_{j\ge3}p_j+\cdots+\suml_{j\ge i}p_j$. Therefore
\beqm
\suml_{i\ge1}p_im_i-m^2_1=\suml_{i\ge2}p_im_i-(m_1-p_1)m_1=\suml_{i\ge2}p_i(m_i-(i-1)m_1)=
\suml_{i\ge2}p_i\Big(\suml_{j\ge2}p_j+\suml_{j\ge3}p_j+\cdots+\suml_{j\ge i}p_j\Big)=\\=
\suml_{i\ge2}p_i\Big(\suml^i_{j\ge2}(j-1)p_j+(i-1)\suml^r_{j>i}p_j\Big)=
\suml_{i\ge2}(i-1)p^2_i+2\suml_{2\le i<j}(i-1)p_ip_j
\end{multline}
Altogether we get:
\beq
\De(Spec)=
(1-\al)^2\frac{\suml_{i\ge2}(i-1)p^2_i+2\suml_{2\le i<j}(i-1)p_ip_j}{2}+(1-\al)\frac{\suml_{i\ge2}(2-i)p_i}{2}+
Frac^{(2)}+Frac^{(1)}_{m_1,\dots,m_r}
\eeq
This gives the difference:
\beq
\De(\mu)-\frac{2}{(1-\al)^2}\De(Spec)=
1-\sum^r_{i=1}p_i+\frac{\suml_{i\ge2}(i-2)p_i}{1-\al}-Frac^{(2)}-Frac^{(1)}_{m_1,\dots,m_r}
\eeq
Finally, remark that $Frac^{(2)}\le\suml_{i\ge1}p_i+r$, while $Frac^{(1)}\le0$. This gives the bound.
\epr


\begin{thebibliography}{99}
\bibitem[AGLV]{AGLV} V.I.Arnol'd, V.V.Goryunov, O.V.Lyashko, V.A.Vasil'ev, {\em Singularity theory. I.}
Reprint of the original English edition from the series Encyclopaedia of Mathematical Sciences
[ Dynamical systems. VI, Encyclopaedia Math. Sci., 6, Springer, Berlin, 1993]. Springer-Verlag, Berlin, 1998. iv+245 pp. ISBN: 3-540-63711-7

\bibitem[Buchweitz-Greuel-1980]{Buchweitz-Greuel-1980} R.-O.Buchweitz, G.-M.Greuel, {\em
The Milnor number and deformations of complex curve singularities.} Invent. Math. 58 (1980), no. 3, 241-�281


\bibitem[Durfee1978]{Durfee1978} A.H. Durfee, {\it The signature of smoothings of complex surface singularities.}
Math. Ann. 232 (1978), no. 1, 85--98.

\bibitem[Givental1983]{Givental} A.B.Givental, {\em The maximum number of singular points on a projective hypersurface.} (Russian) Funktsional. Anal. i Prilozhen. 17 (1983), no. 3, 73–-74


\bibitem[GLS]{GLS} G.M.Greuel, C.Lossen, E.Shustin, {\em Introduction to singularities and deformations.} Springer Monographs in Mathematics. Springer, Berlin, 2007

\bibitem[Kaenders1996]{Kaenders1996} R.Kaenders, {\em The Seifert form of a plane curve singularity determines its intersection multiplicities.}  Indag. Math. (N.S.) 7 (1996), no. 2, 185--197

\bibitem[Kerner.Nemethi.2014]{Kerner.Nemethi.2014} D.Kerner, A.Nemethi, {\em Durfee-type bound for some non-degenerate complete intersection singularities}, in preparation.

\bibitem[Saito1981]{Saito1981}M. Saito, {\em On the exponents and the geometric genus of an isolated hypersurface
 singularity}. Singularities, Part 2 (Arcata, Calif., 1981), 465--472,
 Proc. Sympos. Pure Math., 40, Amer. Math. Soc., Providence, RI, 1983.


\bibitem[Schr.Stee.Stev.]{Schr.Stee.Stev.} R.Schrauwen, J.Steenbrink, J.Stevens, {\em Spectral pairs and the topology of curve singularities.} Complex geometry and Lie theory (Sundance, UT, 1989), 305-�328, Proc. Sympos. Pure Math., 53, Amer. Math. Soc., Providence, RI, 1991.

\bibitem[Steenbrink1976]{Steenbrink1976} J.H.M.Steenbrink, {\em Mixed Hodge structure on the vanishing cohomology.} Real and complex singularities (Proc. Ninth Nordic Summer School/NAVF Sympos. Math., Oslo, 1976), pp. 525--563. Sijthoff and Noordhoff, Alphen aan den Rijn, 1977


\bibitem[Th\`{a}nh-Steenbrink-89]{Thanh-Steenbrink-89}L\^{e} V\v{a}n Th\`{a}nh, J.H.M.Steenbrink, {\em
Le spectre d'une singularit\'{e}� d'un germe de courbe plane.} Acta Math. Vietnam. 14 (1989), no. 1, 87--94

\end{thebibliography}
\end{document}